\documentclass[11pt,letterpaper,reqno]{amsart}
\usepackage{url}
\usepackage{amsthm}
\usepackage{amsmath}
\usepackage{amssymb}
\usepackage{amsfonts}

\usepackage{enumitem}
\usepackage{mathscinet}
\usepackage{lipsum}
\usepackage[english]{babel}
\usepackage[autostyle]{csquotes}
\usepackage{bbold}

\usepackage[english]{babel}

\usepackage[colorlinks=true,
    linkcolor=red,
    citecolor=blue]{hyperref}

\usepackage{graphicx}

\usepackage[utf8]{inputenc}
\usepackage[english]{babel}

\usepackage{color}
\usepackage{comment}

\newtheorem{thm}{Theorem}[section]
\newtheorem{definition}[thm]{Definition}
\newtheorem{theorem}[thm]{Theorem}
\newtheorem*{oseledec*}{Oseledec Theorem}
\newtheorem{cor}[thm]{Corollary}

\newtheorem{lemma}[thm]{Lemma}
\newtheorem{prop}[thm]{Proposition}
\theoremstyle{definition}
\newtheorem{remark}[thm]{Remark}

\makeatletter
\def\moverlay{\mathpalette\mov@rlay}
\def\mov@rlay#1#2{\leavevmode\vtop{%
   \baselineskip\z@skip \lineskiplimit-\maxdimen
   \ialign{\hfil$\m@th#1##$\hfil\cr#2\crcr}}}
\newcommand{\charfusion}[3][\mathord]{
    #1{\ifx#1\mathop\vphantom{#2}\fi
        \mathpalette\mov@rlay{#2\cr#3}
      }
    \ifx#1\mathop\expandafter\displaylimits\fi}
\makeatother

\let\ul\underline
\let\wh\widehat
\let\wt\widetilde

\newcommand{\nocontentsline}[3]{}
\newcommand{\tocless}[2]{\bgroup\let\addcontentsline=\nocontentsline#1{#2}\egroup}

\usepackage{microtype}

\usepackage{wasysym}

\def\Jac{\ensuremath{\mathrm{Jac}}}

\def\H{\ensuremath{\mathrm{H}}}
\title[Thermodynamic Formalism Out of Equilibrium]{Thermodynamic Formalism Out of Equilibrium Part I: Pressure Out of Equilibrium, and
Azuma Inequality and 
Gibbs Processes}

\begin{document}

\author{S. Ben Ovadia, A. Brown, F. Rodriguez-hertz}
\date{}

\begin{abstract}
We study the thermodynamic formalism of topological \\ Markov shifts where the potential varies by a random walk on an exterior system. We define the {\em pressure out of equilibrium} which is associated to such a family of potentials, and prove a corresponding variational principle. We present an application to random dynamical systems. In particular, we study an open condition for random dynamical systems where the randomness is driven by a Gibbs process, and prove hyperbolicity estimates via an Azuma inequality which were previously only known in the i.i.d. setting.
\end{abstract}

\newcommand{\Addresses}{{
  \bigskip
  \footnotesize

  S.~Ben Ovadia, \textsc{Einstein Institute of Mathematics, The Hebrew University of Jerusalem, 91904 Jerusalem, Israel}. \\ \textit{E-mail address}: \texttt{Snir.BenOvadia@mail.huji.ac.il}
  
  A.~Brown, \textsc{Department of Mathematics, Northwestern University,  Evanston, IL 60208, United States}. \\ \textit{E-mail address}: \texttt{awb@northwestern.edu}
  
F.~Rodriguez-Hertz, \textsc{Department of Mathematics, Eberly College of Science, The Pennsylvania State University, State College, PA 16801, United States}. \\
\textit{E-mail address}: \texttt{fjr11@psu.edu}
}}

\maketitle

\tableofcontents

\section{Introduction}

Given an ensemble of particles and a finite collection of states, the Helmholtz free energy is the difference between the potential energy of a distribution, and its entropy scaled by the inverse temperature. Minimizing the Helmholtz free energy is equivalent to maximizing the difference between the entropy and the potential energy scaled by the temperature. In Ruelle's thermodynamic formalism (\cite{Ruelle67}), the variational principle states that given a potential $-\frac{1}{\tau} \phi\in \mathrm{H\ddot{o}l}(\Sigma)$ ($\tau>0$ is the temperature), where $\Sigma$ is a two-sided compact and topologically transitive topological Markov shift,
\begin{equation*}\label{RuellesReduction}
	\max\Big\{h_\nu(T)+\int \phi d\nu: \nu\text{ is }T\text{-inv. and erg.}\Big\}=P(\phi),
\end{equation*}
where $T:\Sigma\to \Sigma$ is the left-shift, and 
\begin{equation}\label{PressDef}
	P(\phi):=\limsup_{n\to\infty}\frac{1}{n}\log \sum_{|\underline{w}|=n
    }e^{\sum_{k=0}^{n-1}\phi\circ T^k(\theta_{\underline{w}})},
\end{equation}
where $\theta_{\underline{w}}\in [\underline{w}]$ maximizes $\sum_{k=0}^{n-1}\phi\circ T^k(\omega)$. Indeed, the value of $P(\phi)$ is 
independent of the choice of $\theta_{\underline{w}}\in [\underline{w}]$
. That is, Ruelle showed that minimizing the free Helmholtz energy over invariant measures (that is over steady states) can be achieved, and that the minimal value can be expressed as the logarithm of the reciprocal of the spectral radius of an associated operator (called the Ruelle operator).

The variational principle is an optimization problem over a collection of measures, with deep physical motivation. Since the pioneering work of Ruelle, several important extensions have been studied, with different related optimization problems. An important extension is in the study of extensions.

In \cite{LedWalt} Ledrappier and Walters study the following optimization. Let  $X$ and $Y$ be two compact metric spaces, let $T:X\to X$ and $S:Y\to Y$ be two continuous maps, and let $\pi:X\to Y$ be a surjection s.t. $\widehat\pi\circ T=S\circ \widehat\pi$. Then, given an $S$-invariant $\nu$ and a potential $\phi\in C(X)$, Ledrappier and Walters optimize $$\sup\Big\{h_\mu(T|S)+\int \phi d\mu: \mu\text{ is }T\text{-invariant and }\mu\circ \widehat\pi^{-1}=\nu \Big\},$$
where $h_\mu(T|S)$ is the relative entropy. That is, they optimize the contribution of the potential together with the entropy coming from the fibers alone, under the constraint that the invariant measures project to a fixed measure $\nu$ on the base.

In \cite{KiferDenkerStadlbauer}, Denker, Kifer, and Stadlbauer study a notion of random topological pressure, and prove that it satisfies a relative variational principle as well (see also \cite{KiferRandomPressure} for a definition of random pressure). In their setting, the fibers are countable Markov shifts, and they study the associated composition of random Ruelle operators w.r.t. to a probability preserving transformation on the base dynamics (see \cite{BG_RandomRuelle} for results in the compact-shift fibers setting). This is a somewhat opposite approach to the one we are taking, as we study randomness driven by a random walk on a fiber which is a general compact metric space. In particular, the notion of pressure which we associate to the random walk on the fiber is not given by a random composition of Ruelle operators by the random walk. In addition, we do not fix a probability measure on the base (nor fibers), but instead seek an optimization over all invariant probabilities. A full description of our setting appears below.

Another extension of Ruelle's variational principle is studied in \cite{BuzziNonlinearTDF}, where the authors study an exponential growth rate as in \eqref{PressDef}, and $\sum_{k=0}^{n-1}\phi\circ T^k$ is generalized to $n\cdot G(\frac{1}{n}\sum_{k=0}^{n-1}\phi\circ T^k)$ with $G$ being a continuous function.

In this paper, we study another extension of Ruelle's variational principle, with applications to random dynamics (in particular to Gibbs processes).

We consider a compact metric space $X$, and a H\"older continuous skew-product over $\Sigma$: $\widehat{F}:\Sigma\times X\to\Sigma\times X$, $\widehat{F}(\omega, t)=(T\omega, F_\omega(t))$, together with an equi-H\"older continuous family of potentials $\phi_t:\Sigma\to \mathbb{R}$, which depend in a H\"older continuous way on $t\in X$. Write $\widehat{\phi}(\omega,t):=\phi_t(\omega)$. We study an optimization problem (i.e. a variational principle) where we allow the potential to be random- depending on the $(X,F_\omega)$-random walk.

We introduce the corresponding {\em pressure out of equilibrium} (POE for short):
\begin{equation}\label{limsupEqIntro}
    \widehat{P}(\widehat{\phi}):=\limsup_{n\to\infty}\sup_{t\in X}\frac{1}{n}\log\sum_{|\underline{w}|=n
    }e^{\sum_{k=0}^{n-1}\phi_{F^k_{\theta_{\underline{w}}}(t)}\circ T^k(\theta_{\underline{w}})},
\end{equation}
where $F^k_\omega:=F_{T^{k-1}\omega}\circ\cdots \circ F_\omega$, and 
$\theta_{\underline{w}}\in [\underline{w}]$ maximizes $\sum_{k=0}^{n-1}\phi_{F^k_{\omega}(t)}\circ T^k(\omega)$. In fact, the limit superior in \eqref{limsupEqIntro} is a limit; And in addition, we sum over all words of length $n$ in our compact setting, however for a possible extensions to non-compact setting, we may restrict to summing over words which begin or end with a fixed letter (see Lemma \ref{limsupIsLim}). Indeed, when $\phi_t(\omega)=\phi(\omega)$, and so the potential does not depend on the $(X,F_\omega)$-random walk, the POE $\wh P(\wh\phi)$ coincides with the usual pressure $P(\phi)$.

 Our POE variational principle (Theorem \ref{VarPrincLowBound}) asserts that when assuming $F_\omega=F_{(\omega_i)_{i\leq0}}$ and 
  {\em $\widehat{\phi}$-holonomies} (which can be verified by an open condition, see Definition \ref{phiHolonomies} and Remark \ref{phiHolonomiesRmk}), we have
\begin{align}\label{mainResult1}
    \widehat{P}(\widehat{\phi})=\max\Big\{&
    \widehat{h}_{\widehat{\nu}}(\widehat{F})+\int \widehat{\phi}d\widehat\nu:
    \widehat\nu\text{ is an 
    }\widehat{F}\text{-inv. prob. on }\Sigma\times X\Big\},
\end{align}
where $\widehat{h}_{\widehat{\nu}}(\widehat{F})$ is the {\em metric entropy out of equilibrium} (see \textsection \ref{EntOutOfEq} and Remark \ref{importRmk}). Some special cases of measures maximizing the POE variational principle are $u$-Gibbs measures and generalized $u$-Gibbs measures (see Remark \ref{gUgibbs}).

One can think of the entropy out of equilibrium as entropy along an unstable direction (from the base dynamics). Moreover, we show that it is dominated by the metric entropy of the projection of the measure to the base dynamics, and the two quantities coincide if and only if the respective measure is {\em future independent}  (Theorem \ref{boundedness}). The future independence property is defined in Definition \ref{theMeas}.

In addition, we prove that when $F_\omega=F_{\omega_0}$, any measure $\tau\in \mathbb{P}(X)$ is $\nu_0$-stationary if and only if $\tau$ is the projection to $X$ of a POE maximizer in \eqref{mainResult1} with $\wh\phi(\omega,t)=\sum_a\log \nu_0(a)\cdot\mathbb{1}_{[a]\times X}$. This means that measures which maximize the POE offer a potential extension of the notion of stationary measures for Gibbs processes rather than just Bernoulli i.i.d. processes, and allow for $F_\omega=F_{(\omega_i)_{i\leq0}}$. See Theorem \ref{stationarity} and Remark \ref{newStationaryRmk}.

We show that $\wh P_{\wh\nu}(\wh \phi)\leq P_{\wh\nu}(\wh \phi)$ for all $\wh F$-invariant measures, and equality holds if and only if $\wh\nu$ is future independent and admits no fiber entropy. In particular, we always have $\wh P(\wh \phi)\leq P(\wh \phi)$. See Theorem \ref{boundedness} and Remark \ref{followsBoundedness}.

When $\wh\phi(\omega,t)=\phi(\omega)$, we show that every POE maximizer is future independent (see Theorem \ref{maxIsFutInd}). In contrast, in \textsection \ref{appeB} we construct an example where no POE maximizer is future independent, and in fact there is a unique maximizer. 

See Remark \ref{remarks2} and Lemma \ref{Sinai} on how to arrange that $F_\omega=F_{(\omega_i)_{i\leq0}}$. Moreover, without the additional assumption of $\widehat{\phi}$-holonomies and $F_\omega=F_{(\omega_i)_{i\leq0}}$, we still establish the ``$\leq$" in \eqref{mainResult1} (Theorem \ref{varPrince}), with the entropy out equilibrium being replaced by the metric entropy of $\widehat{\nu}\circ\pi^{-1}$.

The notion of pressure out of equilibrium arises naturally in the study of skew-products (where we are removing the fiber entropy from the optimization process). The terminology is motivated by the heuristic where the potential on $\Sigma$ keeps varying depending on the random walk on $X$. This leads to a notion of pressure on $\Sigma$ when the potential is ``out of equilibrium". 

The pressure out of equilibrium also arises naturally when studying random dynamical systems. In the recent years significant progress has been made in the study of smooth dynamical systems, where the diffeomorphism is composed randomly. See for example \cite{BenoistQuint,AaronFederico,DeWittDolgopyat}. 

The motivating model where the dynamics (and the potential) depend randomly on an exterior system is an important model, as every physical system in fact interacts with a larger system (hence the laws of motion themselves may depend randomly on a larger system). In this setting, it is natural to assume that the typical behavior of the larger system is given by a physical measure (e.g. SRB) or a measure of maximal entropy- both of which are equilibrium states and are Gibbs measures. Consequently, we are interested in random dynamical systems where the randomness is driven by a Gibbs process. However, almost all previous results treat the i.i.d. case, where the random composition is taken over a finite set of diffeomorphisms.

An important open class of random dynamical systems is the class of systems which satisfy the {\em uniform expansion on average} condition: $\forall \omega^-$,
\begin{equation}\label{uea1st}
    \inf_{(x,\xi)\in T^1M}\int \log |d_xf_{\omega}\xi|d\mu_{\omega^-}(\omega)\geq\chi>0,
\end{equation}
where $f_{\omega}\in \mathrm{Diff}^{1+\alpha}(M)$ ($\omega\mapsto f_{\omega}$ may not be constant on partition sets), $M$ is a closed Riemannian manifold, $\mu$ is a Gibbs measure on a two-sided shift of a H\"older continuous potential $\psi$, 
$\omega^{-}:=(\omega_i)_{i \leq -1}$, 
and $\mu_{\omega^-}$ is the  conditional measure of $\mu$ on the unstable leaf associated with $\omega^-$ (i.e. $\{\omega\in \Sigma: \forall i\leq -1, \omega_i=\omega_i^-\}$, not local in $[\omega_0]:=\{
\omega'=(\omega'_i)_{i\in\mathbb{Z}}\in \Sigma: \omega'_0=\omega_0\}$, but 
the push-forward of the local leaf in $[\omega_{-1}]$). The prevalence of similar conditions has been studied for example in \cite{DimaKrikorian,Rafael,ELS,BEFRH,BergerCarrasco}.

In that setting, we prove the following hyperbolicity estimates (Theorem \ref{final}): $\exists 
C>0$ s.t.
\begin{equation}\label{nonIIDAzuma}
	\forall n\geq 0,\text{ }\sup_{(x,\xi)\in T^1M}\int |d_xf_\omega^n\xi|^{-\beta}d\mu(\omega)\leq Ce^{nP(\psi)-\beta (\chi-O(\beta)) n},
\end{equation}
where $f_\omega^n:=f_{T^{n-1}\omega}\circ \cdots \circ f_\omega$.  

In addition, we prove that \eqref{uea1st} implies an iterated version of \eqref{uea1st} (Corollary \ref{rmkForStrongMixLinBoubd}): $\forall \omega^-$,
\begin{equation}\label{iteratedUEA}
   \inf_{(x,\xi)\in T^1M}\int \log|d_xf_\omega^n\xi|d\mu_{\omega^-}(\omega)\geq nP(\psi)+(\chi-o_n(1)) n.
\end{equation}

One can think of the estimate in \eqref{nonIIDAzuma} as a ``non i.i.d. version of the Azuma inequality" (often also referred to as Hoeffding's lemma), as we provide also a more general form: $\exists C>0$ s.t.
\begin{equation}\label{nonIIDAzuma2}
	\forall n\geq 0,\text{ }\sup_{t\in X}\mathbb{E}\Big(e^{-\beta\sum_{k=0}^{n-1}\varphi_{F_\omega^k(t)}(T^k\omega)}\Big)\leq Ce^{nP(\psi)-\beta (\chi-O(\beta)) n},
\end{equation}
and \eqref{uea1st} translates into
\begin{equation}\label{uea2nd}
   \ \inf_{s\in X}\mathbb{E}(\varphi_s|(\omega_i)_{i\leq -1})\geq \chi>0\text{ a.s.},
\end{equation}
where $\{\varphi_s\}_{s\in X}$ is a family of H\"older continuous potentials $\varphi_s:\Sigma\to \mathbb{R}$ which depend in a H\"older continuous way on $s\in X$,  and $\omega\sim \mu$.  We also prove an iterated form of \eqref{uea2nd}, similarly to \eqref{iteratedUEA} (see Corollary \ref{rmkForStrongMixLinBoubd}). The sequence $\omega$ is sampled w.r.t. $\mu$, and $\{\omega_k\}_{k\geq0}$ is neither i.i.d. nor a Markov chain. We note that applying the Markov inequality to \eqref{nonIIDAzuma2} yields exponential tail estimates for $\sum_{k=0}^{n-1}\varphi_{F_\omega^k(t)}(T^k\omega)$. When $\mu$ is a Gibbs measure, sampling the sequence $\omega$ w.r.t. to $\mu$ is called a {\em Gibbs process}. 

In this setting, even a Markov process is a non-trivial extension of the i.i.d. case (i.e. Bernoulli process). We note that an essential challenge here is that the map $\omega\mapsto F_\omega^n(x)$ is not uniformly H\"older over all $n$, and so standard techniques which apply to martingales do not apply here, even for non-stationary systems (see for example \cite{YeorDima,GouezelChazottes}). An additional key challenge is the fact that already Markov measures (let alone Gibbs measures) no longer enjoy the product-structure that Bernoulli measures enjoy.

\section{Pressure out of equilibrium}\label{POE}

\noindent\textbf{Setup:}
\begin{enumerate}
	\item Let $X$ be a compact metric space.
	\item Let $\Sigma$ be a two-sided compact topological Markov shift, endowed with the left-shift $T:\Sigma\to \Sigma$.
	\item Let $F:\Sigma\to C_\alpha(X,X)$ be a H\"older continuous map, denoted by $\omega\mapsto F_\omega$, where $F_\omega$ is a $\alpha$-H\"older homeomorphism of $X$.
	\item  Given $\widehat{\phi}\in \mathrm{H\ddot{o}l}(\Sigma\times X)$, we denote $\phi_t(\omega)=\widehat{\phi}(\omega,t)$. We view the family $\{\phi_t\}_{t\in X}\subseteq \mathrm{H\ddot{o}l}(\Sigma)$ as a equi-H\"older family of potentials:\\ $\sup_{t\in X}\|\phi_t\|_{\mathrm{H\ddot{o}l}}<\infty$, and $\phi_t$ depends in a H\"older continuous manner on $t\in X$.
 	\item We denote $F_\omega^k:=F_{T^{k-1}\omega}\circ \cdots \circ F_\omega$ for $k\geq1$, and $F_\omega^0=\mathrm{Id}$.
\end{enumerate}

\begin{definition}[Pressure out of equilibrium]\label{DefOfPOE}
	We define the {\em pressure out of equilibrium}, or {\em POE} for short, as $$\widehat{P}(\widehat{\phi}):=\limsup_{n\to\infty}\frac{1}{n}\log\sup_{t\in X}\widehat{Z}_n(\widehat{\phi},t,a),$$ 
	where
	$$\widehat{Z}_n(\widehat{\phi},t,a):=\sum_{|\underline{w}|=n,w_{n-1}=a}e^{\sum_{k=0}^{n-1}\phi_{F^k_{\theta_{\underline{w}}}(t)}(T^k(\theta_{\underline{w}}))},$$
	$[a]\subset \Sigma$, and $\theta_{\underline{w}}\in[\underline{w}]$ maximize $\omega\mapsto \sum_{k=0}^{n-1}\phi_{F^k_{\omega}(t)}(T^k(\omega))$.
\end{definition}

\begin{remark}
    \text{ }
\begin{enumerate}
\item Note, since $\Sigma$ is compact and topologically transitive, the definition of the POE does not depend on $a$. 
Under an assumption of $\widehat{\phi}$-holonomies (see Definition \ref{phiHolonomies}), the definition of the POE also does not depend on the choice of $\theta_{\underline{w}}\in[\underline{w}]$. The condition of $\widehat{\phi}$-holonomies can be satisfied via an open condition of bunching for many interesting examples.
\item Note, $\widehat{Z}_n(\widehat{\phi},t,a)$ is not of the composite form $L_{\phi_n}\circ \cdots \circ L_{\phi_1}\mathbb{1}_{[a]}$! For each word of length $n$, the sequence of potentials in its corresponding weight depends on the word.	
\end{enumerate}
\end{remark}

\begin{lemma}\label{limsupIsLim} The following limit exists:
	$$\widehat{P}(\widehat{\phi})=\lim_{n\to\infty}\frac{1}{n}\log\sup_{t\in X}\widehat{Z}_n(\widehat{\phi},t,a).
   $$
\end{lemma}
\begin{proof}
First, since $\Sigma$ is compact and topologically transitive, there exists $N\in\mathbb{N}$ s.t. for all $[a],[b]\subseteq \Sigma$, for all $n\in\mathbb{N}$,
\begin{align*}
    \widehat{Z}_n(\widehat{\phi},t,a)=&\sum_{\overset{|\underline{w}|=n,}{w_{n-1}=a}}e^{\sum_{k=0}^{n-1}\phi_{ F^k_{\theta_{\underline{w}}}(t)}(T^k(\theta_{\underline{w}}))}\\
\leq& e^{ N\|\widehat{\phi}\|_\infty }\sum_{\overset{|\underline{w}|=n+N,}{w_{n+N-1}=b}}e^{\sum_{k=0}^{n+N-1}\phi_{ F^k_{\theta_{\underline{w}}}(t)}(T^k(\theta_{\underline{w}}))}\\
=&\widehat{Z}_{n+N}(\widehat{\phi},t,b)\cdot e^{N\|\widehat{\phi}\|_\infty}.
\end{align*}
Then, for all $[a']\subseteq \Sigma$,
\begin{align}\label{Feb19}
    \sup_t \widehat{Z}_n(\widehat{\phi},t,a')\leq& \sup_t\sum_{[a]\subseteq\Sigma} \widehat{Z}_{n}(\widehat{\phi},t,a)\\
\leq&e^{N\|\widehat{\phi}\|_\infty}\sup_t\sum_{[a]\subseteq\Sigma} \widehat{Z}_{n+N}(\widehat{\phi},t,a')\nonumber\\
\leq& \#\{[a]\subseteq \Sigma\}\cdot e^{N\|\widehat{\phi}\|_\infty} \sup_t \widehat{Z}_{n+N}(\widehat{\phi},t,a').\nonumber
\end{align}
Then, 

$$\frac{1}{n}\log\sup_t \widehat{Z}_n(\widehat{\phi},t,a)\text{ converges} \Leftrightarrow\frac{1}{n}\log \widehat{Z}_n(\widehat{\phi})\text{ converges,}$$ 
where $$\widehat{Z}_n(\widehat{\phi}):= \sup_t\sum_{[a]\subseteq\Sigma} \widehat{Z}_{n}(\widehat{\phi},t,a)=\sup_t \sum_{|\underline{w}|=n}e^{\sum_{k=0}^{n-1}\phi_{F^k_{\theta_{\underline{w}}}(t)}(T^k(\theta_{\underline{w}}))} .$$ Moreover, the two sequences share their set of limits points.
	
Then we notice that for all $n,m\in \mathbb{N}$, $$\log \widehat{Z}_{n+m}(\widehat{\phi})\leq \log \widehat{Z}_n(\widehat{\phi})+ \log \widehat{Z}_m(\widehat{\phi}),$$
since the collection of admissible words of length $n+m$ is smaller or equal to the concatenation of all words of length $n$ and all words of length $m$, and $\sup_t \{A_t\cdot B_t\}\leq \sup_t A_t\cdot \sup_t B_t $. Then by Fekete's lemma, we are done. 
\end{proof}

\section{Entropy out of Equilibrium}\label{EntOutOfEq}

\begin{definition}
     $\Sigma^-:=\{(\omega_i)_{i\leq-1 }:\omega\in\Sigma\}$. Elements of $\Sigma^-$ are often denoted by $\omega^-$. Given $\omega^-\in \Sigma$, its associated {\em unstable leaf} is $$V^u_{\leq-1}(\omega^-):=\{\omega\in \Sigma: \forall i\leq -1, \ \omega_i=\omega^-_i\}.$$
\end{definition}

\begin{definition}
   We define the following partitions:
   \begin{enumerate}
    \item $\mathcal{C}:=\{[a]: [a]\subseteq \Sigma\}$,
       \item $\widehat{\mathcal{C}}:=\{[a]\times X: [a]\subseteq \Sigma\}$,
       \item For $q\in \mathbb{N}$, $\mathcal{C}_q:=\bigvee_{i=0}^{q-1}T^{-i}\mathcal{C}$,
        \item For $q\in \mathbb{N}$, $\widehat{\mathcal{C}}_q:=\bigvee_{i=0}^{q-1}\widehat{F}^{-i}\widehat{\mathcal{C}}$.
   \end{enumerate}
\end{definition}

\begin{remark}
    Note, $\widehat{F}$ preserves the fiber structure, and so elements of $\widehat{C}_q$ are saturated by fibers.
\end{remark}

\begin{definition}[Metric entropy out of equilibrium]
    Given an $\widehat{F}$-invariant Borel probability $\widehat{\eta}$ on $\Sigma^-\times X$, we define its {\em metric entropy out of equilibrium} by 
    \begin{equation}\label{EOEform}
        \widehat{h}_{\widehat{\eta}}(\widehat{F}):=\limsup_{n\to\infty} \int \frac{1}{n}\H_{\eta_{(\omega^-,t)}}\Big(\widehat{\mathcal{C}}_q\Big)d\tau(\omega^-,t),
    \end{equation}
    where $\widehat{\eta}=\int \eta_{(\omega^-,t)} d\tau(\omega^-, t)$ is the disintegration of $\widehat{\eta}$ into conditional measures by the measurable partition $\{V^u_{\leq -1}(\omega^-)\times\{t\}\}_{\omega^-\in \Sigma^-, t\in X}$.
\end{definition}

\begin{definition}[Metric pressure out of equilibrium]
    Given an $\widehat{F}$-invariant Borel probability $\widehat{\eta}$ on $\Sigma\times X$, and a  potential $\widehat{\phi}:\Sigma\times X\to\mathbb{R}$, we define its {\em metric pressure out of equilibrium} by 
    $$\widehat{P}_{\widehat{\eta}}(\widehat{\phi}):=\widehat{h}_{\widehat{\eta}}(\widehat{F})+\int \widehat{\phi}d\widehat{\eta}.$$
\end{definition}

\begin{definition}[Local entropy out of equilibrium]
   Given an $\widehat{F}$-invariant Borel probability $\widehat{\eta}$, the {\em local entropy out of equilibrium} of $\widehat{\eta}$ is defined for $\tau$-a.e. $(\omega,t)$ as,
   \begin{equation}\label{locEntForm}\widehat{h}_{\widehat{\eta}}((\omega,t)):=\limsup_{n\to\infty}\frac{-1}{n}\log\eta_{(\omega^-,t)}([\omega]_n),
   \end{equation}
   where $[\omega]_n:=[\omega_0,\ldots,\omega_{n-1}]$,  $\omega^-=(\omega_i)_{i\leq -1}$, and  $\widehat{\eta}=\int \eta_{(\omega^-,t)} d\tau(t)$ is the disintegration of $\widehat{\eta}$ into conditional measures by the measurable partition $\{V^u_{\leq -1}(\omega^-)\times\{t\}\}_{\omega^-\in \Sigma^-, t\in X}$.
\end{definition}

\begin{lemma}\label{strongUEnt}
    Assume that $F_\omega=F_{(\omega_i)_{i\leq0}}$, and let $\widehat{\eta}$ be an ergodic $\widehat{F}$-invariant Borel probability. Then the limit superior in \eqref{EOEform} is in fact a limit, and for $\widehat{\eta}$-a.e. $(\omega,t)$,
    $$
        \widehat{h}_{\widehat{\eta}}(\widehat{F})=\lim_{n\to\infty}\frac{-1}{n}\log\eta_{(\omega^-,t)}([\omega]_n),
    $$
     $\widehat{\eta}=\int \eta_{(\omega^-,t)} d\tau(t)$ is the disintegration of $\widehat{\eta}$ into conditional measures by the measurable partition $\{V^u_{\leq -1}(\omega^-)\times\{t\}\}_{\omega^-\in \Sigma^-, t\in X}$.
\end{lemma}
\begin{proof}
    When $F_\omega=F_{(\omega_i)_{i\leq0}}$, the partition $\xi^u:=\{V^u_{\leq -1}(\omega^-)\times\{t\}\}_{\omega^-\in \Sigma^-, t\in X}$ is an increasing partition: 
    $$ \widehat{F}^{-1}\xi^u>\xi^u\text{ and }\bigvee_{n\geq0}\widehat{F}^{-n}\xi^u=\text{partition into points}.$$
    Indeed, for $(\omega,t)\in \Sigma\times X$ and $n\in\mathbb{N}$, we have
    $$\xi^u((\omega,t))\cap[\omega]_n=\xi^u_n((\omega,t)),$$
where $\xi_n^u:=\widehat{F}^{-n}\xi^u$.

    Therefore, $\widehat{h}_{\widehat{\eta}}((\omega,t))$ is precisely the conditional entropy of $\widehat{F}$ along the measurable and increasing partition of $\xi^u$. In this case, the limit superior in \eqref{locEntForm} is in fact a limit (by the point-wise ergodic theorem, and for integrability note that $\H_{\eta_{(\omega^-,t)}}(\widehat{F}^{-1}\xi^u)\leq \log\#\mathcal{C}$), and it is an invariant function and thus constant almost everywhere. It then follows easily from the dominated convergence theorem that for $\widehat{\eta}$-a.e. $(\omega',s)$,
    \begin{align*}
    \widehat{h}_{\widehat{\eta}}(\widehat{F})=&\lim_n\int \frac{1}{n}\H_{\eta_{(\omega^-,t)}}(\xi^u_n)d\tau(\omega^-,t)\\
    =&\int\lim_n \frac{1}{n}\H_{\eta_{(\omega^-,t)}}(\xi^u_n)d\tau(\omega^-,t)=\int\widehat{h}_{\widehat{\eta}}(\omega,t) d\widehat{\eta}=\widehat{h}_{\widehat{\eta}}(\omega',s).
        \end{align*}
\end{proof}

\begin{remark}\label{importRmk} Under the assumptions of Lemma \ref{strongUEnt} (i.e. $F_\omega=F_{(\omega_i)_{i\leq0}}$), the entropy out of equilibrium is in fact the conditional entropy corresponding the increasing measurable partition of strong unstable leaves (i.e. the unstable direction of the symbolic part).
\end{remark}

\begin{definition}[Entropy out of equilibrium]
We define the {\em entropy out of equilibrium}:    $$\wh h (\wh F):=\wh P(0).$$
\end{definition}

\begin{remark}\label{rmkEOE}
    In the setting of Theorem \ref{VarPrincLowBound}, we get that 
    \begin{equation*}
        \wh h (\wh F)=\max\Big\{\wh h_{\wh\eta}(\wh F): \wh\eta\text{ is }\wh F\text{-inv.} \Big\}.
    \end{equation*}
By Proposition \ref{maxIsFutInd}, we get in addition that, $$\wh h(\wh F)=h_\mathrm{top}(T).$$  
\end{remark}

\section{Future Independence, Stationarity, and Entropy Out of Equilibrium}\label{FutInd}
\begin{definition}[Future independence]\label{theMeas}
    Let $\widehat{\nu}$ be a Borel probability measure. Write $$\widehat{\nu}=\int \widehat{\nu}_{\omega^-} d\nu^-(\omega^-),$$
a disintegration by conditional probabilities w.r.t. the partition $\{V^u_{\leq-1}(\omega)\times X\}_{\omega\in\Sigma}$, where $V^u_{\leq-1}(\omega)=\{\omega'\in \Sigma: \omega'_i=\omega_i,\forall i\leq-1\}$, and $\nu^-$ is the projection of $\nu$ to $\Sigma^-$.

We say that  $\widehat{\nu}$ is  {\em future independent} if for $\nu^-$-a.e. $\omega^-$, 
$$\widehat{\nu}_{\omega^-}=\sigma_{\omega^-}\times \nu_{\omega^-},$$
where $\nu=\int \nu_{\omega^-}d\nu^-$ is the disintegration of $\nu$ by symbolic one-step-longer local unstable leaves $V^u_{\leq-1}(\cdot)$, and $\sigma_{\omega_{\leq -1}}$ is a Borel probability measure on $X$.
\end{definition}

\begin{remark}\label{remark1}\text{ }
\begin{enumerate}
\item     
    We will frequently impose the standing condition that the generator of our cocycle $\omega\mapsto F_\omega$ is constant along sets of the form $V^u(\omega)$.  That is, we assume 
$$F_\omega=F_{(\omega_i)_{i\leq 0}}.$$
The generators of the inverse cocycle $F^{-1}$ are  $\omega\mapsto (F_{T^{-1}\omega})^{-1}$ are then constant on sets of the form $$V^u_{\leq-1}(\omega) = \{ \omega'\in \Sigma: \omega'_i=\omega_i \  \text {for all $i\le -1$}\}. $$
That is, if $\omega \mapsto G_\omega = (F_{T^{-1}\omega})^{-1}$ generates the inverse cocycle, then  $G_\omega$ is of the form $$G_\omega =G_{(\omega_i)_{i\leq -1}} .$$
    \item We note that the collection of measures which are future independent can be identified with the collection of $\widehat{G}$-invariant measures, where $\widehat{G}:\Sigma^-\times X\to\Sigma^-\times X $ and $\Sigma^{-}:=\{(\omega_i)_{i\leq-1}:\omega\in \Sigma\}$.
    
    Also, note that a measure is future independent if and only if its fiber conditional measures are invariant under holonomies by unstable strong leaves.

\item An i.i.d.\ random walk on a space $X$ and a stationary measure $\nu_0$ on $X$  naturally induce a skew-extension over a full shift $\Sigma$ with an invariant measure $\widehat\nu$.   In this case, the  base measure $\mu=\pi_*\widehat\nu$ is always a Bernoulli measure.  

The stationarity of $\nu_0$ induces geometric properties of the disintegration $\widehat\nu_\omega$.  The family of measures $\omega\mapsto\widehat\nu_\omega$ is constant along local unstable manifolds (i.e. elements of $\Sigma^-$), hence independent of the future, as in the choice of terminology ``future independence".
\item The future independence property can also be described as the product-structure of unstable-fiber conditional measures.
\item The same notion can be viewed geometrically as a type of {\em stationarity}: The conditional measures on the fiber $X$ are stationary when varying the stable leaf. Indeed, when $\nu$ is a Bernoulli measure, and $f_\omega=f_{\omega_0}$, the standard notion of a $\nu$-stationary measure (i.e. $\sigma =\int (f_\omega)_*\sigma d\nu(\omega)$ where $\sigma$ is the projection of $\widehat{\nu}$ onto the fiber) coincides with the notion of future independence. See Theorem \ref{stationarity}.
\item Notice that if $f_\omega=f_{\omega_0}$, the condition on $\sigma$, $\sigma =\int (f_\omega)_*\sigma d\nu(\omega)$,
only takes into account the $\nu$-measure of cylinders of length $1$ and treats $\nu$ as a Bernoulli measure. Thus, it does not truly extend the notion of stationarity well for Gibbs processes.

\end{enumerate}
 \end{remark}

\begin{theorem}\label{boundedness}
    Given an $\widehat{F}$-invariant Borel probability $\widehat{\eta}$ on $\Sigma\times X$, 
    $$\widehat{h}_{\widehat{\eta}}(\widehat{F})\leq h_\eta(T),$$
    where $\eta:=\widehat{\eta}\circ \pi^{-1}$ and $\pi:\Sigma\times X\to\Sigma$ is the projection. Moreover, equality holds if and only if $\wh \eta$ is future independent.
\end{theorem}
\begin{proof}
    Let $\wh\eta=\int \wh\eta_{\omega^-}d\eta^-$ be the disintegration of $\wh\eta$ w.r.t. the measurable partition $\{\{\omega^-\}\times X:\omega^-\in\Sigma^-\}$. Given $\wh\eta_{\omega^-}$ on $V^u_{\leq -1}(\omega^-)\times X$, disintegrate further into $\wh\eta_{\omega^-}=\int \eta_{(\omega^-,t)}d\tau_{\omega^-}(t)$. Write $\Phi(x)=-x\log x$. Then $\Phi$ is concave, and so by Jensen's inequality, 
    \begin{align*}
        \H_{\wh\eta_{\omega^-}}\Big(\wh{\mathcal{C}}
        \Big)=&\sum_{A\in \wh{\mathcal{C}}
        }\Phi(\wh\eta_{\omega^-}(A))=
\sum_{A\in \wh{\mathcal{C}}
}\Phi(\int\eta_{(\omega^-,t)}(A)d\tau_{\omega^-}(t))\\
\geq&\sum_{A\in \wh{\mathcal{C}}
}\int\Phi(\eta_{(\omega^-,t)}(A))d\tau_{\omega^-}(t)=\int \sum_{A\in \wh{\mathcal{C}}
}\Phi(\eta_{(\omega^-,t)}(A))d\tau_{\omega^-}(t)\\
=&\int \H_{\eta_{(\omega^-,t)}}\Big(\wh{\mathcal{C}}
\Big)d\tau_{\omega^-}(t),
\end{align*}
and equality holds if and only if $\eta_{(\omega^-,t)}(A)$ is constant for $\tau_{\omega^-}$-a.e. $t$, for all $A\in\wh{\mathcal{C}}
$. Since this needs to hold 
for $\eta^-$-a.e. $\omega^-$, otherwise we get a strict drop in entropy: We get that $$\widehat{h}_{\widehat{\eta}}(\widehat{F})=\int\int \H_{\eta_{(\omega^-,t)}}\Big(\wh{\mathcal{C}}
\Big)d\tau_{\omega^-}(t)d\eta^-\leq \int \H_{\wh\eta_{\omega^-}}\Big(\wh{\mathcal{C}}
        \Big)d\eta^- = h_\eta(T),$$ 
and equality holds if and only $\wh\eta_{\omega^-}$ coincide on the $\sigma$-algebra generated by $\wh{\mathcal{C}}$ for $\tau_{\omega^-}$-a.e. $t$, for $\eta^-$-a.e. $\omega^-$. By the invariance of $\widehat{\eta}$, the conditionals must coincide on the $\sigma$-algebra generated by $\bigvee_{i=0}^{\infty}\wh F^{-i}\wh{\mathcal{C}}$, which refines into points in local unstable leaves.
\end{proof}

\begin{remark}\label{followsBoundedness}
    From Theorem \ref{boundedness} the following statement follows: Given an $\wh F$-invariant probability $\wh\eta$, we always have $\wh P_{\wh\eta}(\wh \phi)\leq P_{\wh\eta}(\wh \phi):=h_{\wh \nu}(\wh F)+\int \wh \phi d\wh \eta$ and equality holds if and only if $\wh\nu$ is future independent, and the fiber entropy is zero.
\end{remark}

\section{Variational Principle: Upper Bound}\label{varPcple}

The proof of the variational principle for H\"older continuous potentials on TMSs goes through Ruelle operators. We cannot use the techniques of Ruelle operators in this setting, as the POE is not the spectral radius of a Ruelle operator, nor of a random composition of Ruelle operators given by the random walk.

\begin{theorem}[Upper bound]\label{varPrince}
	$$\widehat{P}(\widehat{\phi})\leq \max\Big\{h_{\widehat{\nu}\circ \pi^{-1}}(T)+\int \widehat{\phi} d\widehat\nu: \widehat\nu\text{ is an }\widehat F\text{-inv. Borel prob.}\Big\},$$
    where $\pi:\Sigma\times X\to\Sigma$ is the projection.
\end{theorem}

This proof borrows ideas from \cite{Walters}.
\begin{proof}
First we note that it is enough to bound $\widehat{P}(\widehat{\phi})$ by a single element in the set on the l.h.s. as by expansiveness the maximum exists. Let $n\in \mathbb{N}$, and let $t_n \in X$ s.t. $\widehat{Z}_n(\widehat{\phi},t_n)=\max_t\widehat{Z}_n(\widehat{\phi},t)$. Set,
\begin{equation}\label{forNuCompt}
\widetilde{\nu}_n:=\sum_{|\underline{w}|=n}\frac{e^{\widehat{\phi}^{(n)}(\theta_{\underline{w}},t_n)}}{\widehat{Z}_n(\widehat{\phi})}\delta_{\theta_{\underline{w}}}.
\end{equation}

Then,
\begin{align}\label{newUpperBoundJensen}
    & \sum_{|\underline{w}|=n}-\log\widetilde{\nu}_n([\underline{w}])\cdot \widetilde{\nu}_n([\underline{w}])\\
    =& \sum_{|\underline{w}|=n}\log\frac{e^{\widehat{\phi}^{(n)}(\theta_{\underline{w}},t_n)}}{\widetilde{\nu}_n([\underline{w}])}\cdot \widetilde{\nu}_n([\underline{w}])-\int \widehat{\phi}^{(n)}(\omega,t_n)d\widetilde{\nu}_n(\omega)
    \nonumber\\
    = &\log \sum_{|\underline{w}|=n}e^{\widehat{\phi}^{(n)}(\theta_{\underline{w}},t_n)}-n\cdot \int \widehat{\phi}d\widehat{\nu}_n\nonumber\\
    \leq &\max_{(\omega^-,t)}\log \widehat{Z}_n(\widehat{\phi},\omega^-,t)-n\cdot \int \widehat{\phi}d\widehat{\nu}_n\nonumber,
\end{align}
where
$$\widehat{\nu}_n:=\frac{1}{n}\sum_{k\leq n-1}(\widetilde{\nu}_n\times\delta_{t_n})\circ \widehat{F}^{-k}.$$

To bound the l.h.s. of \eqref{newUpperBoundJensen} from above, we first notice that it can be written as follows:
\begin{align}\label{AprApr}
   \sum_{|\underline{w}|=n}-\log\widetilde{\nu}_n([\underline{w}])\cdot \widetilde{\nu}_n([\underline{w}])
=\H_{\widetilde{\nu}_n}\Big(\mathcal{C}_n\Big).
\end{align}

We continue to bound from above the r.h.s. of \eqref{AprApr}. Fix $1\leq q < n$, and for any $0\leq j\leq q-1$ set $a(j):=[\frac{n-j}{q}]$. Then, 
 $$\mathcal{C}_{n}=\bigvee_{r=0}^{a(j)-1}T^{-rq+j}\mathcal{C}_{q} \vee \bigvee_{l\in R}T^{-l}\mathcal{C},$$
 where $\# R\leq 2q$. And so,
\begin{align}\label{eq7}
\H_{\widetilde{\nu}_n}\Big(\mathcal{C}_n\Big)\leq&\sum_{r=0}^{a(j)-1}\H_{\widetilde{\nu}_n}\Big(T^{-rq+j}\mathcal{C}_{q}\Big)+ \H_{\widetilde{\nu}_n}\Big(\bigvee_{l\in R}T^{-l} \mathcal{C}\Big)\nonumber\\
	 \leq& \sum_{r=0}^{a(j)-1}\H_{\widetilde{\nu}_n\circ T^{-(rq+j)}}\Big(\mathcal{C}_{q}\Big)+ 2q\log\#\mathcal{C}.
\end{align}
Summing \eqref{eq7} over $0\leq j\leq q-1$ and dividing by $n\cdot q$, we get 
\begin{align}\label{eq8}
\frac{1}{n}\H_{\widetilde{\nu}_n}\Big(\mathcal{C}_n\Big)\leq&\frac{1}{n}\sum_{j=0}^{q-1}\sum_{r=0}^{a(j)-1}\frac{1}{q}\H_{\widetilde{\nu}_n\circ T^{-(rq+j)}}\Big(\mathcal{C}_{q}
\Big)+ \frac{2q\log\#\mathcal{C}}{n}\nonumber\\
\leq&\frac{1}{q}\H_{\frac{1}{n}\sum_{p\leq n}\widetilde{\nu}_n\circ T^{-k}}\Big(\mathcal{C}_{q}
\Big)+ \frac{2q\log\#\mathcal{C}}{n}\nonumber\\
=&\frac{1}{q}\H_{\widehat{\nu}_n\circ \pi^{-1}}\Big(\mathcal{C}_{q}\Big)+ \frac{2q\log\#\mathcal{C}}{n}.
\end{align}

Then, sending $n_i\to\infty$, then taking limit superior on $q\to\infty$, we get
\begin{align}\label{eq9}
    \limsup_{i\to\infty}\frac{1}{n_i}\H_{\widetilde{\nu}_{n_i}}\Big(\mathcal{C}_{n_i}\Big)\leq& \limsup_q\lim_i\frac{1}{q}\H_{\widehat{\nu}_{n_i}\circ \pi^{-1}}\Big(\mathcal{C}_q\Big)\nonumber\\
    =&\limsup_q\frac{1}{q}\H_{\widehat{\nu}\circ \pi^{-1}}\Big(\mathcal{C}_q\Big)=h_{\widehat{\nu}\circ \pi^{-1}}(T),
\end{align}
where $\pi$ is continuous and so $\widehat{\nu}_n\circ\pi^{-1}\to\widehat{\nu}\circ\pi^{-1}$.
\end{proof}

\section{Variational Principle}

\begin{definition}[$\phi$-holonomies]\label{phiHolonomies}
    We say that $(\Sigma\times X,\widehat{F})$ admits {\em $\widehat{\phi}$-holonomies}, if $\exists C>0$ s.t. $$\sup_{n\geq0,t\in X}\sup_{\omega,\omega'\in \Sigma: \omega_i=\omega_i', i\in [0,n]}|\widehat{\phi}^{(n)}(\omega,t)-\widehat{\phi}^{(n)}(\omega',t)|\leq C.$$
\end{definition}

\begin{remark}\label{phiHolonomiesRmk}\text{ }
\begin{enumerate}
\item We will mostly discuss Definition \ref{phiHolonomies} in the context of H\"older continuous potentials, as in our setting.
    \item Stable and unstable holonomies (see Definition \ref{StableHolonomies}) together imply $\widehat{\phi}$-holonomies. One can verify this by adding and subtracting the Birkhoff average over the word $[\omega,\omega']$ (i.e., the Smale bracket). In particular, stable and unstable holonomies are an open condition for a large class of partially hyperbolic systems (see Remark \ref{HolsAreOpen}).

\end{enumerate}
     \end{remark}

\begin{definition}
    Given $\omega^-\in \Sigma^-$ and $t\in X$,
$$\widehat{Z}_n(\widehat{\phi},\omega^-,t):=\sum_{\overset{|\underline{w}|=n,}{\omega^-_{-1}\to w_0}}e^{\widehat{\phi}^{(n)}(\omega^-\cdot \theta_{\underline{w}},t)}.$$
Similarly, for $\omega\in \Sigma$, we define
$$\widehat{Z}_n(\widehat{\phi},\omega,t):=\sum_{\overset{|\underline{w}|=n,}{\omega_{-1}\to w_0}}e^{\widehat{\phi}^{(n)}([\omega, \theta_{\underline{w}}],t)},$$
where here $[\cdot,\cdot]$ denotes the Smale brackets over the $-1$-th position.
\end{definition}

\begin{remark}
    In Theorem \ref{VarPrincLowBound} below we assume $\wh\phi$-holonomies. This assumption is used only in order to compare the following growth rates:
    $$\frac{1}{n}\max_{(\omega^-,t)}\log\widehat{Z}_n(\widehat{\phi},\omega^-,t)\asymp\frac{1}{n}\log \widehat{Z}_n(\widehat{\phi}).$$
 To avoid this assumption, one may study the quantity on the l.h.s.   
\end{remark}

\begin{definition}[Space of measured leaves]\label{MeasuredLeaves}
    We define the following {\em space of measured leaves}:
$$\mathcal{S}:=\Big\{(\omega^-,t,\nu)\in \Sigma^-\times X\times \mathbb{P}(\Sigma\times X):\mathrm{Supp}(\nu)\subseteq V^u_{\leq-1}(\omega^-)\times\{t\}\Big\}.$$
\end{definition}

\begin{remark}\label{MeasuredLeavesRmk}
    Note that $\mathcal{S}$ is compact. This fact becomes important for us in the proof of Theorem \ref{VarPrincLowBound}. 
\end{remark}

\begin{theorem}[POE variational principle]\label{VarPrincLowBound} 
Assume that $(\Sigma\times X,\widehat{F})$ admits $\widehat{\phi}$-holonomies. Assume further that $F_\omega=F_{(\omega_i)_{i\leq 0}}$. 
 Then,
$$\widehat{P}(\widehat{\phi})= \max\Big\{\widehat{h}_{\widehat{\eta}}(\widehat{F})+\int \widehat{\phi} d\widehat\eta: \widehat\eta\text{ is an }\widehat{F}\text{-inv. Borel prob.}\Big\}.$$
\end{theorem}
\begin{proof} \text{ }

\medskip
\textbf{Step 1:} Using the additional assumptions, we prove a tighter upper bound of the POE than Theorem \ref{varPrince}.

\medskip
Let $n\in\mathbb{N}$. Fix $\widetilde{\omega}^{-,n}\in \Sigma^-$ which is not periodic and s.t. 
$$\max_t\sum_{|\underline{w}|=n, \wt\omega^{-,n}_{-1}\to w_0}e^{\wh\phi^{(n)}(\theta_{\underline{w}},t)}\geq \frac{1}{\#\mathcal{C}}\wh Z_n(\wh\phi).$$
Let $\widetilde{t}_n\in X$ which maximizes $t\mapsto \wh Z_n(\wh \phi,\wt \omega^{-,n},t)$. Note, by $\wh\phi$-holonomies we have that $\wh Z_n(\wh \phi,\wt \omega^{-,n},\wt t_n)=e^{\pm C}\cdot \sum_{|\underline{w}|=n, \wt\omega^{-,n}_{-1}\to w_0}e^{\wh\phi^{(n)}(\theta_{\underline{w}},\wt t_n)}$. We define the following (probability) leaf measure which is supported on $V^u_{\leq-1}(\wt \omega^-)\times\{\wt t_n\}$:
\begin{align*}
\widehat{\nu}_{(\wt \omega^{-,n},\wt t_n)}^{(n)}:=&\frac{1}{\wh{Z}_n(\wh\phi,\wt\omega^{-,n},\wt t_n)}\sum_{\wh F^{-n}((\omega^{-},t))=(\wt\omega^{-,n},\wt t_n)}e^{\wh \phi^{(-n)}(\omega^{-}, t)}\delta_{\theta_{\omega^-}}\times\delta_{\wt t_n},
\end{align*}
where $\theta_{\omega^-}\in \widehat{F}^{-n}[V^u_{\leq-1}(\omega^-)]\subseteq V^u_{\leq-1}(\wt\omega^-)$.

We define the following ambient measure as an average over leaf measures:
$$\wh\nu_n:=\frac{1}{n}\sum_{k\leq n}\wt\nu_{(\wt\omega^-,\wt t_n)}^{(n)}\circ \wh F^{-k}.$$

Recall Definition \ref{MeasuredLeaves} and Remark \ref{MeasuredLeavesRmk}, and note that we can write 
$$\widehat{\nu}_n=\int \nu d\Xi_n((\omega^-,t,\nu)),$$
where $\Xi_n\in \mathbb{P}(\mathcal{S})$ is a Borel probability measure on the (compact) space of measured leaves.

The map $\Pi:(\omega^-,t,\nu)\mapsto \nu$ is continuous (as a projection), and so if $\Xi_{n_j}\to \Xi$, then $$\widehat{\nu}_{n_j}\to \int \nu d\Xi((\omega^-,t,\nu))\equiv\wh\nu,$$ which is an $\wh F$-invariant probability measure. Moreover, $\nu\mapsto H_\nu\Big(\wh{\mathcal{C}}\Big)$ is continuous as $\wh {\mathcal{C}}$ is a clopen partition, and so
$$\int \H_\nu\Big(\wh{\mathcal{C}}\Big)d\Xi_{n_j}((\omega^-,t,\nu))\to \int\H_\nu\Big(\wh{ \mathcal{C}}\Big)d\Xi((\omega^-,t,\nu)).$$

\medskip
\textbf{Claim:} 
\begin{equation}\label{propertyN}
    \int \H_\nu\Big(\wh{\mathcal{C}}\Big) d\Xi_n((\omega^-,t,\nu))\geq -\int\wh\phi d\wh\nu_n+\frac{1}{n}\log\wh Z_n(\wh\phi,\wt\omega^{-,n},\wt t_n).
\end{equation}

\medskip
\textbf{Proof:} Since $\wt\omega^{-,n}$ is not periodic, for every $(\omega^-,t,\nu)$ in the support of $\Xi_n$ there exists a unique $k=k_\nu\leq n$ and $(\omega^{-}_\nu,t_\nu)$ s.t. $\widehat{F}^{-k}((\omega^{-}_\nu,t_\nu))=(\wt\omega^{-,n},\wt t_n)$ and s.t. $$\nu\propto\sum_{\widehat{F}^{-(n-k)}((\omega^{-},t))=(\omega^{-}_\nu,t_\nu)}e^{\wh\phi^{(-(n-k))}(\omega^-,t)}\delta_{\wh F^k(\theta_{\omega^-},\wt t_n)}.$$
Indeed, the normalizing quantity for $\nu$ must be $\wh Z_{n-k}(\wh\phi,\omega^{-,'},t')$. Then, let $\nu\in \mathrm{Supp}(\Xi_n\circ\Pi)$ and let $k=k_\nu$, and so
\begin{align*}
\H_{\nu}\Big(\wh{\mathcal{C}}\Big)=&\sum_{\overset{\wh F((\omega^-,t))}{=(\omega^{-}_\nu,t_\nu)}}-\frac{e^{\wh\phi(\omega^-_\nu,t_\nu)}\wh Z_{n-k-1}(\wh\phi,\omega^-,t)}{\wh Z_{n-k}(\wh\phi,\omega^-_\nu,t_\nu)}\log\frac{e^{\wh\phi(\omega^-_\nu,t_\nu)}\wh Z_{n-k-1}(\wh\phi,\omega^-,t)}{\wh Z_{n-k}(\wh\phi,\omega^-_\nu,t_\nu)}\\
=&\int -\wh\phi d\nu-\int\log\frac{\wh Z_{n-k-1}(\wh\phi,\omega,t)}{\wh Z_{n-k}(\wh\phi,\omega^-_\nu,t_\nu)}d\nu\\
=&\int -\wh\phi d\nu+\int\log\frac{\wh Z_{n-k}\circ \wh F^k(\wh\phi,T^{-k}\omega,\wt t_n)}{\wh Z_{n-k-1}\circ \wh F^{k+1}(\wh\phi,T^{-k}\omega,\wt t_n)}d\nu.
\end{align*}

Integrating both sides by $\Xi_n$, we get
\begin{align}\label{strongerClaim}
     \int \H_\nu\Big(\widehat{\mathcal{C}}\Big)d\Xi_n((\omega^-,t,\nu))=&\int -\widehat{\phi}d\wh\nu_n\nonumber\\
     +&\frac{1}{n}\int \log \prod_{k=0}^{n-1}\frac{\wh Z_{n-k}\circ \wh F^k(\wh\phi,\omega,\wt t_n)}{\wh Z_{n-k-1}\circ \wh F^{k+1}(\wh\phi,\omega,\wt t_n)}d\wt\nu_{(\wt\omega^{-,n},\wt t_n)}^{(n)}\nonumber\\
     =&\int -\widehat{\phi}d\wh\nu_n+\frac{1}{n}\log\wh Z_n(\wh\phi,\wt\omega^{-,n},\wt t_n)\nonumber\\
      \geq & \int -\widehat{\phi}d\wh\nu_n+\frac{1}{n}\log\wh Z_n(\wh\phi)-\frac{\log C+\log \#\wh{\mathcal{C}}}{n}.
\end{align}
This concludes the claim.

\medskip
From \eqref{strongerClaim} we get that,
$$\int \H_\nu\Big(\wh{\mathcal{C}}\Big)d\Xi((\omega^-,t,\nu))\geq \int \wh\phi d\wh\nu+\wh P(\wh\phi).$$
Let $\Gamma:\Sigma^-\times X\times \mathbb{P}(\Sigma\times X)\to \Sigma^-\times X$ be the projection map. Write $\xi:=\Xi\circ \Gamma^{-1}$. Then, $\wh \nu=\int \nu_{(\omega^-,t)} d\xi((\omega^-,t))$, where for $\xi$-a.e. $(\omega^-,t)$,
$$\nu_{(\omega^-,t)}=\int \nu d\Xi_{(\omega^-,t)}(\nu),$$
where
$$\Xi=\int \Xi_{(\omega^-,t)} d\xi((\omega^-,t)).$$
Now notice, for $\xi$-a.e. $(\omega^-,t)$, 
\begin{align*}
   \H_{\int \nu d\Xi_{(\omega^-,t)}(\nu)}\Big(\wh{\mathcal{C}}\Big)\geq& \int \H_{\nu}\Big(\wh{\mathcal{C}}\Big)d\Xi_{(\omega^-,t)}(\nu).
\end{align*}
Integrating both sides by $\xi$, we get
\begin{align}\label{infoFunchhat}
    \int \H_{\nu_{(\omega^-,t)}}\Big(\wh {\mathcal{C}}\Big)d\wh\nu=&\int \H_{\int \nu d\Xi_{(\omega^-,t)}(\nu)}\Big(\wh{\mathcal{C}}\Big)d\xi((\omega^-,t))\\
    \geq&\int \int \H_{\nu}\Big(\wh{\mathcal{C}}\Big)d\Xi_{(\omega^-,t)}(\nu)d\xi((\omega^-,t))\nonumber\\
    =&\int \H_{\nu}\Big(\wh{\mathcal{C}}\Big)d\Xi\geq-\int \wh\phi d\wh\nu+\wh P(\wh\phi).\nonumber
\end{align}
We now notice that the l.h.s. of \eqref{infoFunchhat} equals $\int \wh I_{\wh \nu}(\omega,t)d\wh\nu$, where $\wh I$ is the information function of the conditional entropy on unstable leaves of the base dynamics (recall Theorem \ref{strongUEnt}). Therefore, $$\wh P_{\wh \nu}(\wh \phi)\geq \wh P(\wh \phi).$$

\medskip
\textbf{Step 2:} We continue to bound the POE from below by the metric POE for all $\widehat{F}$-invariant probability measures. In this case,
$$\widehat{P}_{\widehat{\nu}}(\widehat{\phi})\geq\widehat{P}(\widehat{\phi})\geq \sup\{\widehat{P}_{\widehat{\eta}}(\widehat{\phi})\}_{\widehat{\eta}},$$
and so the supremum is in fact a maximum.

\medskip
\textbf{Proof:} Let $\widehat{\eta}$ be an $\widehat{F}$-invariant Borel probability measure. Let $\widehat{\eta}=\int \eta_{(\omega^-,t)}d\tau(\omega^-,t)$ be the disintegration of $\widehat{\eta}$ into conditional measures by the measurable partition $\{V^u_{\leq-1}(\omega^-)\times\{t\}\}_{\omega^-\in\Sigma^-,t\in X}$.

\begin{align}\label{newUpperBoundJensen2}
    &\int \sum_{\overset{|\underline{w}|=n,}{\omega^-_{-1}\to w_0}}-\log\eta_{(\omega^-,t)}([\underline{w}])\cdot \eta_{(\omega^-,t)}([\underline{w}])d\tau(\omega^-,t)\\
    \leq&\int \Big(\sum_{\overset{|\underline{w}|=n,}{\omega^-_{-1}\to w_0}}\log\frac{e^{\widehat{\phi}^{(n)}(\omega_{\underline{w},\omega^-,t},t)}}{\eta_{(\omega^-,t)}([\underline{w}])}\cdot \eta_{(\omega^-,t)}([\underline{w}])-\int \widehat{\phi}^{(n)}(\omega,t)d\eta_{(\omega^-,t)}(\omega)\Big)d\tau\nonumber\\
    \leq &\int \Big(\log \sum_{\overset{|\underline{w}|=n,}{\omega^-_{-1}\to w_0}}e^{\widehat{\phi}^{(n)}(\omega_{\underline{w},\omega^-,t},t)}\Big)d\tau(\omega^-,t)-n\cdot \int \widehat{\phi}d\widehat{\eta}\nonumber\\
    \leq &\max_{(\omega^-,t)}\log \widehat{Z}_n(\widehat{\phi},\omega^-,t)-n\cdot \int \widehat{\phi}d\widehat{\eta}\nonumber,
\end{align}
where $\omega_{\ul w, \omega^-,t}$ maximizes $\wh\phi^{(n)}$ on $(V^u_{\leq -1}(\omega^{-})\cap[\ul w])\times \{t\}$ the and the second inequality is due to Jensen. And so,
$$\int \frac{1}{n}\H_{\eta_{(\omega^-,t)}}\Big(\bigvee^{n-1}\widehat{F}^{-i}\widehat{\mathcal{C}}\Big)d\tau(\omega^-,t)\leq \frac{1}{n}\log \widehat{Z}_n(\widehat{\phi})-\int \widehat{\phi}d\widehat{\eta}.$$
By the definition of the metric entropy out of equilibrium, the limit superior on the l.h.s. equals to $\widehat{h}_{\widehat{\eta}}(\widehat{F})$. By 
Definition \ref{DefOfPOE} the limit superioir on the r.h.s. converges to $\widehat{P}(\widehat{\phi})-\int \widehat{\phi}d\widehat{\eta}$.
\end{proof}

\begin{remark}\label{gUgibbs}\text{ }
\begin{enumerate}
    \item To compare Theorem \ref{VarPrincLowBound} with Theorem \ref{varPrince}, in \textsection \ref{appeB} we construct an example where $$\wh P(\wh \phi)<\max\Big\{h_{\wh\eta\circ \pi^{-1}}(T)+\int \wh\phi d\wh\eta:\wh\eta\text{ is }\wh F\text{-invariant}\Big\}.$$ 
\item A special case of measures which obtain the maximum in the variational principle of Theorem \ref{VarPrincLowBound} are $u$-Gibbs measures, and generalized $u$-Gibbs measures (see \cite{GuGibbs}). Consider two toral Anosov diffeomorphisms $A,B:\mathbb{T}^2\to\mathbb{T}^2$. Let $F_x:\mathbb{T}^2\to\mathbb{T}^2$, $x\in \mathbb{T}^2$, be a family of perturbations of $B$ which preserve a cone family s.t. $x\mapsto F_x$ is constant on local unstable leaves of $A$. Let $\wh F:\mathbb{T}^2\times \mathbb{T}^2\to\mathbb{T}^2\times \mathbb{T}^2$ be defined by $\wh F(x,y)=(A x, F_x(y))$, then $\wh F$ is a uniformly hyperbolic map. Let $\wh\phi(x,y)=-\log \Jac(d_xA|_{E^u(x)})$. We may code the base dynamics and think of the base as a topological Markov shift. In this case, every measure which obtains the maximum in Theorem \ref{VarPrincLowBound} is a generalized $u$-Gibbs (as its conditionals on the unstable leaves corresponding to the base dynamics are smooth). To see that the POE is $0$ in this case, one carries out the construction of the lower bound in the proof of Theorem \ref{VarPrincLowBound} with a smooth leaf measure rather than $\wt \nu^{(n)}_{\cdot}$, as in \cite{uGibbsPesinSinai}. If we assume further that $A$ dominates all the hyperbolic maps on the fiber (by, for example, speeding up the base map), then in this case any maximizer in the POE variational principle is a $u$-Gibbs measure. In both cases the system is uniformly hyperbolic and so it admits an SRB measure, which will be a POE maximizer. Below, in Theorem \ref{maxIsFutInd}, we prove that such measures must be future independent.
\end{enumerate}
\end{remark}

\section{POE Maximizers, and Future Independence and Stationarity}

\begin{theorem}\label{maxIsFutInd}
In the setting of Theorem \ref{VarPrincLowBound}, if $\wh\phi(\omega,t)=\phi(\omega)$ then $\wh P(\wh\phi)=P(\phi)$, and any measure which obtains the maximum in the POE variational principle is future independent.
\end{theorem}
\begin{proof}
If $\wh\phi(\omega,t)=\phi(\omega)$, then
$$\wh P(\wh\phi)=\lim\frac{1}{n}\log \max_t\sum_{|\underline{w}|=n}e^{\wh\phi^{(n)}(\theta_{\underline{w}},t)}=\lim\frac{1}{n}\log \sum_{|\underline{w}|=n}e^{\phi^{(n)}(\theta_{\underline{w}})}=P(\phi).$$

Then, by Theorem \ref{boundedness}, it follows that for any maximizer $\wh \eta$ in the POE variational principle, 
\begin{align*}
    h_\eta(T)+\int \phi d\eta\geq\wh h_{\wh\eta}(\wh F)+\int \phi d\eta=\wh P(\wh \phi)=P(\phi)=\max\Big\{h_{\eta'}(T)+\int \phi d\eta'\Big\}, 
\end{align*}
and so $\wh h_{\wh\eta}(\wh F)=h_\eta(T)$.

Finally, by Theorem \ref{boundedness} once more, it follows that $\wh\eta$ is future independent.
\end{proof}

\begin{remark}
    Note, $\wh P(\wh \phi)=P(\phi)\leq P(\wh\phi)$. The inequality may be strict, depending on whether or not the fiber dynamics admit positive entropy.
\end{remark}

\begin{theorem}[POE maximizers and stationary measures]\label{stationarity}
    Let $\nu_0\in \mathbb{P}(\mathcal{C})$, and write $\psi=\sum_{a}\mathbb{1}_{[a]}\log \nu_0([a])$. Assume that $F_\omega=F_{\omega_0}$. 
    Then,
    \begin{enumerate}
        \item Every $\nu_0$-stationary measure can be obtained as the projection of a POE maximizer to $X$,
        \item Every POE maximizer projects to a $\nu_0$-stationary measure on $X$.
    \end{enumerate}
\end{theorem}
\begin{proof}\text{ }
\begin{enumerate}
    \item  
    Let $\tau$ be a $\nu_0$-stationary measure. Let $\nu^+$ be a $\psi$-equilibrium state (i.e. the associated Bernoulli measure $\nu_0^\mathbb{N}$) on $\Sigma^+:=\{(\omega_i)_{i\geq0}:\omega\in\Sigma\}$. Then, $\wh\nu^+:=\nu^+\times\tau$ is $\wh F^+:\Sigma^+\times X\to \Sigma^+\times X$ invariant
   . Let $\wh\nu$ be the natural extension of $\wh\nu^+$ to an $\wh F$-invariant measure 
   on the two-sided system $\Sigma\times X$.  
   We show that $\wh P_{\wh \nu}(\psi)=0$, and so $\wh\nu$ is a maximizer.

 Indeed, a consequence of the construction of $\widehat \nu$ via Martingale convergence (see \cite{Martingales}) implies that $$\widehat \nu = \lim_{n\to \infty} \widehat{F}^n_*(\nu_0^{\mathbb Z}\times \tau).$$  Each measure $\widehat{F}^n_*(\nu_0^{\mathbb Z}\times \tau)$, while not $\widehat{F}$-invariant, is future independent.  The future independence property is closed among all measures with a fixed marginal $\nu$. Thus $\wh\nu$ is future independent. By Theorem \ref{boundedness}, $\wh h_{\wh \nu}(\wh F)=h_\nu(T)$, and so $\wh P_{\wh\nu}(\psi)=h_{\nu}(T)+\int \psi d\nu^+=0$.
   
   \item 
If $\wh\nu$ is a maximizer in the POE variational principle, then by Theorem \ref{maxIsFutInd} and Theorem \ref{boundedness}, $h_{\nu^+}(T)=
   \wh h_{\wh\nu}(\wh F)$, where $\nu^+$ which is the projection of $\wh\nu$ to $\Sigma^+$. Therefore, $\nu^+$ is a maximizer for the pressure of $\psi$. By the invariance of $\wh\nu^+$ (the projection of $\wh\nu$ to $\Sigma^+\times X$), we get that $\wh\nu\circ \Upsilon^{-1}$ is $\nu_0$-stationary, where $\Upsilon:\Sigma\times X\to X$ is the projection.
\end{enumerate}
\end{proof}

\begin{remark}\label{newStationaryRmk}
    Theorem \ref{stationarity} implies that the terminology ``$\nu_0$-stationary", which is equivalent to ``$\psi$-stationary" when $F_\omega=F_{\omega_0}$, can be extended to H\"older continuous potentials and to allow for $F_\omega=F_{(\omega_i)_{i\leq0}}$. Moreover, Theorem \ref{maxIsFutInd} says that when $\psi\in \mathrm{H\ddot{o}l}(\Sigma)$, then any maximizer in the POE variational principle is future independent. This fact adds structure to the extended notion of stationary measures. In \textsection \ref{appeB} we construct an example where no POE maximizer is future independent, when the potential depends on both $\omega$ and $t$.
\end{remark}

\section{Hyperbolicity estimates for uniformly expanding on average Gibbs processes}\label{expansionSect}

In this section we apply the results of \textsection \ref{varPcple} in the setting of Gibbs processes.

\noindent\textbf{Setup:}
\begin{enumerate}
	\item Let $M$ be a closed Riemannian manifold of dimension $d\geq 2$.
	\item Let $f:\Sigma\to \mathrm{Diff}^{1+\alpha}(M)$ be a H\"older continuous map, denoted by $\omega\mapsto f_\omega$.
\item Let $\mu$ be a $T$-invariant probability measure on $\Sigma$, and assume further that $\mu$ is a Gibbs measure for the potential $\psi\in \mathrm{H\ddot{o}l}(\Sigma)$.
\item Given $(x,\xi)=t\in X$ where $X:= T^1M$, write $\varphi_t(\omega):=\log |d_xf_\omega\xi|$, and note that 
$\sup_{t\in X}\|\varphi_\omega\|_{\mathrm{H\ddot{o}l}}<\infty$.
\item Given $\omega\in \Sigma$, set $V^u_{\leq-1}(\omega)=\{\omega'\in \Sigma:\forall i\leq-1, \omega_i'=\omega_i\}$.
\item Let $\Sigma^-:=\{(\omega_i)_{i\leq-1}: \omega\in \Sigma\}$.  We view 
$\Sigma^-$ as a parametrization of the collection  of local unstable sets $\omega\mapsto V^u_{\leq-1}(\omega)$.  Given $\omega^-\in \Sigma^-$, we write
 $\mu_{\omega^-}$ for the  conditional measure of $\mu$ on the local unstable set associated with $\omega^-$ (i.e., the set $\{\omega\in \Sigma: \forall i\leq -1, \omega_i=\omega_i^-\}$).
    \item Assume w.l.o.g. that 
\begin{enumerate}
    \item $P(\psi)=0$, where $P(\cdot)$ denotes the Gurevich pressure,
    \item $\psi(\omega)=\psi((\omega_i)_{i\leq-1})$,
    \item $L_\psi1=1$, where $L_\psi$ is the associated Ruelle operator on $C(\Sigma^-)$.
\end{enumerate}
This can be arranged by considering the H\"older continuous $(\psi-P(\psi))$-cohomologue $\log\mu_{\omega_{\leq -2}}(V^u_{\leq -1}(\omega))\leq 0$ where $\{\mu_{\omega_{\leq -2}}\}_{\omega\in \Sigma}$ are the conditionals of $\mu$ on two-steps-longer local unstable leaves.	
\end{enumerate}

\begin{definition}[Uniform expansion on average]\label{UEACOND}
    We say that $\mu$ is {\em uniformly $\widehat{\varphi}$-expanding on average} if:
    \begin{equation*}
	    \exists\chi>0\text{ s.t. }\forall \omega^-\in \Sigma^-, \ \inf_{t\in X}\int\varphi_\omega(t)d\mu_{\omega^-}(\omega)\geq \chi.
	\end{equation*}
\end{definition}

\begin{remark}
\text{ }
\begin{enumerate}
\item We do not require $\omega\mapsto f_\omega$ to be constant on partition sets.
\item The random dynamical system $f_\omega\sim \mu$ is called a {\em Gibbs process}. It naturally extends Bernoulli processes and Markov processes.
\end{enumerate} 
\end{remark}

Gibbs processes translate naturally to the setting of thermodynamic formalism out of equilibrium (recall \textsection \ref{POE}):
$$F_\omega(x,\xi):=\Big(f_\omega(x),\frac{d_xf_\omega\xi}{| d_xf_\omega\xi|}\Big),$$
where $X=T^1M$.
 
\medskip
\noindent We wish to prove that there exist $\beta,\gamma>0$ s.t. for all $n$ sufficiently large, for all $(x,\xi)\in T^1M$,
 \begin{equation}\label{eq1}
 	\int_\Sigma |d_xf_\omega^n\xi|^{-\beta}d\mu(\omega)\leq e^{-\gamma n}.
 \end{equation}

\begin{remark}
    
 By the Markov inequality, \eqref{eq1} implies that for all $(x,\xi)\in T^1M$, every word $\omega\in \Sigma$ expands uniformly, aside for an exponentially small exceptional set of words. Such estimates have proved themselves to be the fundamental estimates which are used in the study of statistical properties of random dynamical system (see for example \cite{DimaKrikorian,DeWittDolgopyat}). Previously, such estimates were shown for i.i.d. dynamical systems (i.e. $\mu$ is a Bernoulli measure, and $\omega\mapsto f_\omega$ is constant on partition sets). However, the physical motivation suggests that the randomness which drives our dynamics should be given by a Gibbs measure, hence motivating our extension.
\end{remark}

\subsection{Reduction to thermodynamic formalism out of equilibrium}\label{rightHere}

\begin{theorem}\label{thm2}
Assume that $(\Sigma,T)$ is topologically transitive, and that $(\Sigma\times T^1M,\widehat{F})$ admits $\widehat{\varphi}$-holonomies where $\widehat{\varphi}(\omega,(x,\xi)):=\log|d_xf_\omega\xi|$. Then, for all $\beta\in (0,1)$, for all $n\in \mathbb{N}$ and all $t=(x,\xi)\in T^1M$,
	$$\int |d_xf_{\omega}^{n}\xi|^{-\beta}d\mu= \widehat{C}_{\psi,\varphi}^{\pm1}\cdot \widehat{Z}_n(\widehat{\phi},t,a),$$
	where, $\phi_t:=\psi-\beta\varphi_t$, $\widehat{C}_{\psi,\varphi}\geq1$ is a global constant and $[a]\subset \Sigma$.
\end{theorem}
\begin{proof}
Let $n\in\mathbb{N}$ and $t=(x,\xi)\in T^1M$. For all $[a]\subset\Sigma$ fix some $\omega_a\in [a]$, and let $\theta_{\underline{w}}\in [\underline{w}]$ be some element of the cylinder which maximizes $\sum_{k=0}^{n-1}\phi_{F^k_{\omega}(t)}(T^k(\omega))$. Then,
\begin{align}\label{eq2}
 		\int |d_xf_\omega^{n}\xi|^{-\beta}d\mu=&\int e^{-\beta \widehat{\varphi}^{(n)}(\omega,t)}d\mu=C_\varphi^{\pm 1}\sum_{|\underline{w}|=n}\mu([\underline{w}])e^{-\beta \widehat{\varphi}^{(n)}(\theta_{\underline{w}},t)}
        \\
 		=& (C_\psi C_\varphi)^{\pm 1}\sum_{|\underline{w}|=n}e^{\psi^{(n)}(\theta_{\underline{w}})-\beta \widehat{\varphi}^{(n)}(\theta_{\underline{w}},t)}\nonumber\\
        =& (C_\psi C_\varphi)^{\pm 1}\sum_{|\underline{w}|=n}e^{\widehat{\phi}^{(n)}(\theta_{\underline{w}},t)}\nonumber\\
 		= &( C_\psi C_\varphi)^{\pm 1}\sum_{[a]\subset \Sigma}C_a^{\pm1}\sum_{|\underline{w}|=n ,w_{n-1}=a}e^{\widehat{\phi}^{(n)}(\underline{w}\cdot \omega_a,t)},\nonumber
\end{align}
where $\omega_a\in [a]$ is any element of $[a]$, and $C_a:=N_a\cdot \sup_t \|\phi_t\|_\infty$ where $\forall [b]\subset \Sigma$, $a\xrightarrow[]{n_{ab}}b \xrightarrow[]{n_{ba}}a$ with $n_{ab}, n_{ba}\leq N_a$.

\end{proof}

\begin{remark}
    \text{ }
\begin{enumerate}
	\item The proof Theorem \ref{thm2} in fact implies the following more general statement:
$$\forall t\in X, \int e^{-\beta\sum_{k=0}^{n-1}\varphi_{F^k_\omega(t)}(T^k(\omega))} d\mu_\psi(\omega)= \widehat{C}_{\psi,\varphi}^{\pm1}\cdot\widehat{Z}_n(\psi-\beta\widehat{\varphi},t,a).$$
\item An uppr-bound alone in \eqref{eq2} with no lower-bound, holds without the assumption of $\widehat{\varphi}$-holonomies.
\item Theorem \ref{thm2} demonstrates that our definition of the pressure out of equilibrium (recall \textsection \ref{POE}) is in fact optimal. If one wishes to gain hyperbolicity estimates such as in \eqref{eq1}, then the POE is the optimal quantity to study.
\end{enumerate}

\end{remark}

\begin{cor}\label{cor1} Under the assumptions of Theorem \ref{thm2},
$$\limsup\frac{1}{n}\log \sup_t \int e^{-\beta\widehat{\varphi}^{(n)}(\omega,t)}d\mu(\omega)=\widehat{P}(\widehat{\phi}).$$
and in particular, when $\widehat{\varphi}(\omega,t)=\log |d_{x_t}f_\omega\xi_t|$ and $X=T^1M$,
$$\limsup\frac{1}{n}\log \sup_t \int |d_xf_{\omega}^{n}\xi|^{-\beta}d\mu\leq \widehat{P}(\widehat{\phi}).$$
\end{cor}

\begin{remark}
    By Corollary \ref{cor1}, proving that the POE is negative for all $\beta>0$ sufficiently small, would prove \eqref{eq1}. 
\end{remark}

\subsection{Bounding the POE}\label{POEleqTP}

\begin{definition}[Stable and unstable holonomies]\label{StableHolonomies}
  Given $\omega,\omega'\in \Sigma$, let $I_{\omega\omega'}:\{\omega\}\times X\to \{\omega'\}\times X$ be the natural identification of the fibers. We say that $(\Sigma\times X,\widehat{F})$ admits {\em stable holonomies} if 
  $\forall \omega,\omega'\in\Sigma\text{ s.t. } \omega_i=\omega'_i,\forall i\geq 0$, the limit $$\lim_{k\to\infty}\widehat{F}^{-k}\circ I_{T^k\omega T^k\omega'}\circ \widehat{F}^k\in \mathrm{H\ddot{o}l(X:X)}
    $$
    exists exponentially fast uniformly in $\omega,\omega'$. Similarly,  we say that $(\Sigma\times X,\widehat{F})$ admits {\em unstable holonomies} if $\forall \omega,\omega'\in\Sigma\text{ s.t. } \omega_i=\omega'_i,\forall i\leq 0$, the limit
    $$\lim_{k\to\infty}\widehat{F}^{k}\circ I_{T^{-k}\omega T^{-k}\omega'}\circ \widehat{F}^{-k}\in \mathrm{H\ddot{o}l(X:X)}
    $$
    exists exponentially fast uniformly in $\omega,\omega'$.
\end{definition}

\begin{remark}\label{HolsAreOpen}
    Partially hyperbolic systems which are given by a perturbation of a fibered system over a uniformly hyperbolic systems admit both stable and unstable holonomies under an open condition of bunching. In turn, this implies $\widehat{\phi}$-holonomies as well (see Definition \ref{phiHolonomies}).
\end{remark}

\begin{definition}[Weak unstable holonomies]\label{WeakHolonomies}
    We say that $(\Sigma\times X,\widehat{F})$ admits {\em weak unstable holonomies} if $$\sup_{t\in X}\sup_{\omega_i=\omega'_i,\forall i\leq -1}d(F_{T^{-j}\omega}^j(t),F_{T^{-j}\omega'}^j(t))=o_j(1),$$
    where a sequence $\{a_j\}_{j\geq0}$ is denoted as $o_j(1)$ if it tends to $0$ as $j$ tends to $\infty$.
\end{definition}

\begin{remark}\label{remarks2} In particular, the condition of unstable holonomies implies that one can conjugate $\widehat{F}$ to a system $\widehat{F}^-$ which admits weak unstable holonomies (see \textsection \ref{appe}). 
\end{remark}

\begin{theorem}\label{final}
    Assume that $(\Sigma\times X,\widehat{F})$ admits weak unstable holonomies and 
    $\widehat{\varphi}$-holonomies, and that $\mu$ is uniformly $\widehat{\varphi}$-expanding on average (as in Definition \ref{UEACOND}) with $\chi>0$. Assume also that $(\Sigma,T)$ is topologically transitive. Then, 
    $$\widehat{P}(\widehat{\phi})\leq -\beta (\chi- O(\beta)).$$
\end{theorem}
We remind the reader that $\widehat{\phi}:=\psi-\beta\widehat{\varphi}$, and so the l.h.s. depends on $\beta$ as well.
\begin{proof}
For every $T$-inv. Borel probability $\nu$, $P_\nu(\psi)\leq P(\psi)=0$ (recall \textsection \ref{expansionSect}). Then, by Theorem \ref{varPrince}, 
    $$\widehat{P}(\widehat{\phi})\leq-\beta \int \widehat{\varphi}d\widehat{\nu}=-\beta \lim_{k\to\infty}\int \frac{1}{n_k}\widehat{\varphi}^{(n_k)}(\omega,t_{n_k}) d\widetilde{\nu}_{n_k}(\omega),$$
 where $\widetilde{\nu}_n$ is given by \eqref{forNuCompt}, and $n_k\uparrow\infty$.

 We continue to show that $\lim_{k\to\infty}\int \frac{1}{n_k}\widehat{\varphi}^{(n_k)}(\omega,t_{n_k}) d\widetilde{\nu}_{n_k}(\omega)\geq \chi - O(\beta)$, by showing that for all $n$, $\int \frac{1}{n}\widehat{\varphi}^{(n)}_{t_{n}} d\widetilde{\nu}_n\geq \chi- O(\beta)$.

The first step of the proof is to provide a substitute to the measure $\widetilde{\nu}_n$ with which it is slightly easier to work, but which does not change the final statement. We provide a substitute for the measure $\widetilde{\nu}_n:$ $$\widetilde{\nu}_n:=\frac{1}{\sum\limits_{|\underline{w}|=n}e^{-\beta\widehat{\varphi}^{(n)}(
\theta_{\underline{w}},t_n)}\mu([\underline{w}])}\sum_{|\underline{w}|=n}\delta_{\theta_{\underline{w}}} e^{-\beta\widehat{\varphi}^{(n)}(
\theta_{\underline{w}},t_n)}\mu([\underline{w}]).$$ 
We first observe that the new form is still a probability measure. Indeed, by the Gibbs estimates for $\mu$, $\mu([\underline{w}])=e^{\psi^{(n)}(\theta_{\underline{w}})}C_\psi^{\pm1}$. Therefore, since a global multiplicative factor does not change the value of $\limsup\frac{1}{n}\log Z_n(\widehat{\phi})$, we may use the new form of $\widetilde{\nu}_n$ in the proof of Theorem \ref{varPrince}. It is useful, as it provides an ad-hoc form in which one can view $\widetilde{\nu}_n$ as given by a kernel over $\mu$, with some additional nice properties.

Let $K_n$ be the $\mu$-kernel defined by $$K_n(\omega):=\frac{e^{-\beta \widehat{\varphi}^{(n)}(\theta_{(\omega_{-n+1},\ldots,\omega_{0})},t_n)}}{\sum\limits_{|\underline{w}|=n}e^{-\beta\widehat{\varphi}^{(n)}(\theta_{\underline{w}},t_n)}\mu([\underline{w}])}.$$
One can check that indeed $\int K_n d\mu=1$. Moreover, $K_n$ is constant on sets of the form $[w_{-n+1},\ldots, w_0]$. And finally, for all such cylinder $[\ul w]=[w_{-n+1},\ldots, w_0]$, we have $(K_n\cdot\mu)([\ul w])=\wt \nu_n([\ul w])$. Therefore, for all $n\in \mathbb{N}$,
\begin{equation}\label{sourceOf1OverN}
    \int \frac{1}{n}\wh \varphi^{(n)}(\cdot,t_n) d(K_n \mu)=\int \frac{1}{n}\wh \varphi^{(n)}(\cdot,t_n) d\wt\nu_n\pm\frac{C_{\wh\varphi}}{n}.
\end{equation}

For any $j\leq n$, we continue to bound  $\int\widehat{\varphi}\circ \widehat{F}^j(\omega,t_{n})d(K_n \cdot\mu)$: Notice that $K_n$ is constant up to a factor of $e^{\pm \beta(\|\widehat{\varphi}\|_\infty+C_{\widehat{\varphi}})}$ on sets of the form $\{\omega'\in \Sigma: \widetilde{\omega}_i=\omega_i',\forall i\leq -1\}$, $\widetilde{\omega}\in \Sigma$ (i.e. one-step-longer local unstable leaves). Then, we write $K_n(\omega)=K_n(\omega_{\leq-1})\cdot (1+\Delta_\beta(\omega))$ where $|\Delta_\beta|\leq2\beta(\|\widehat{\varphi}\|_\infty+C_{\widehat{\varphi}})$ for all $\beta>0$ small enough, and so
\begin{align}\label{forRemarkAfter}
& \int\widehat{\varphi}\circ \widehat{F}^j(\omega,t_{n})  K_n(\omega)d\mu(\omega)\\
    =& \int\widehat{\varphi}(T^j\omega,F^j_{\omega}(t_{n})) K_n(\omega)d\mu(\omega)\nonumber\\
    =& \int \int\widehat{\varphi}(T^j\omega,F^j_{\omega}(t_{n})) K_n(\omega)d\mu_{\widetilde{\omega}_{\leq-1}}(\omega)d\mu(\widetilde{\omega}_{\leq-1})\nonumber\\
=& \int \int K_n(\widetilde{\omega}_{\leq-1})(1+\Delta_\beta(\omega))\widehat{\varphi}(T^j\omega,F^j_{\omega}(t_{n})) d\mu_{\widetilde{\omega}_{\leq-1}}(\omega)d\mu(\widetilde{\omega}_{\leq-1})\nonumber\\
    =& \int K_n(\widetilde{\omega}_{\leq-1})\int\widehat{\varphi}(T^j\omega,F^j_{\omega}(t_{n})) d\mu_{\widetilde{\omega}_{\leq-1}}(\omega)d\mu(\widetilde{\omega}_{\leq-1})\nonumber\\
+&\int K_n(\widetilde{\omega}_{\leq-1})\int\Delta_\beta(\omega)\cdot \widehat{\varphi}(T^j\omega,F^j_{\omega}(t_{n})) d\mu_{\widetilde{\omega}_{\leq-1}}(\omega)d\mu(\widetilde{\omega}_{\leq-1})\nonumber.
\end{align}

We bound the bottom summand in the last line, denoted by $\star$:
\begin{align*}
  |\star|\leq &2\beta(\|\widehat{\varphi}\|_\infty+C_{\widehat{\varphi}})\|\widehat{\varphi}\|_\infty\int K_n(\widetilde{\omega}_{\leq-1}) d\mu(\widetilde{\omega}_{\leq -1})\\
=&2\beta(\|\widehat{\varphi}\|_\infty+C_{\widehat{\varphi}})\|\widehat{\varphi}\|_\infty e^{\pm \beta(\|\widehat{\varphi}\|_\infty+C_{\widehat{\varphi}}) }\int K_n d\mu\leq 3\beta(\|\widehat{\varphi}\|_\infty+C_{\widehat{\varphi}})\|\widehat{\varphi}\|_\infty,
\end{align*}
for all $\beta>0$ small enough w.r.t. $\widehat{\varphi}$.

We continue to bound $\int\widehat{\varphi}(T^j\omega,F^j_{\omega}(t_{n})) d\mu_{\widetilde{\omega}_{\geq-1}}(\omega)$ from below for every $\widetilde{\omega}_{\leq-1}$. By the conformality of the conditional measures,
\begin{align}\label{strongMixLinBound}
  &\int\widehat{\varphi}(T^j\omega,F_\omega^j(t_n)) d\mu_{\widetilde{\omega}_{\leq-1}}(\omega)\nonumber\\
  =& 
\sum_{T^j_R\tau_{\leq-1}=\widetilde{\omega}_{\leq -1}}e^{\psi^{(j)}
  (\tau_{\leq -1})} \int \widehat{\varphi}(\omega,F_{T^{-j}\omega}^j(t_n)) d\mu_{\tau_{\leq -1}}(\omega)\nonumber\\
  =& 
\sum_{T^j_R\tau_{\leq-1}=\widetilde{\omega}_{\leq -1}}e^{\psi^{(j)}
  (\tau_{\leq -1})} \Big(\int \widehat{\varphi}(\omega,s_{j,\tau_{\leq-1},n}) d\mu_{\tau_{\leq -1}}(\omega)-o_j(1)\Big)\nonumber\\
\geq&\sum_{T^j_R\tau_{\leq-1}=\widetilde{\omega}_{\leq -1}}e^{\psi^{(j)}(\tau_{\leq -1})} (\chi-o_j(1))\nonumber\\
=&(L_\psi^j1)(\widetilde{\omega}_{\leq -1})\cdot (\chi-o_j(1))=\chi-o_j(1),
\end{align}
where $T_R$ is the right-shift, $s_{j,\tau_{\leq-1},n}:=F_{T^{-j}\omega_\tau}^{j}(t_n)$, $\omega_\tau\in \mathrm{Supp}(\mu_{\tau_{\leq -1}})$ is constant up to a $o_j(1)$ on the support of $\mu_{\tau_{\leq-1}}$ by the weak unstable holonomies, and the inequality is by the uniform $\widehat{\phi}$-expansion on average condition. Plugging this back in \eqref{forRemarkAfter} yields (together with the observation that $\int K_n(\widetilde{\omega}_{\leq-1})d\mu =e^{\pm\beta(\|\widehat{\varphi}\|_\infty+C_{\widehat{\varphi}}) }$),
\begin{align}\label{Lastly}
    &\int\widehat{\varphi}\circ \widehat{F}^j(\omega,t_{n})  K_n(\omega)d\mu(\omega)\\\geq& e^{-\beta(\|\widehat{\varphi}\|_\infty+C_{\widehat{\varphi}})}(\chi-o_j(1))-3\beta(\|\widehat{\varphi}\|_\infty+C_{\widehat{\varphi}})\|\widehat{\varphi}\|_\infty,\nonumber
\end{align} for all $j\leq n$. Recalling \eqref{sourceOf1OverN}, averaging \eqref{Lastly} over all $j\leq n$ and letting $n_k\to\infty$, we get that  $\int \widehat{\varphi}d\widehat{\nu}\geq e^{-\beta(\|\widehat{\varphi}\|_\infty+C_{\widehat{\varphi}})}\chi-3\beta(\|\widehat{\varphi}\|_\infty+C_{\widehat{\varphi}})\|\widehat{\varphi}\|_\infty=\chi-O(\beta)$.
 
\end{proof}

\begin{remark}\text{ }
\begin{enumerate}
    \item One can view Theorem \ref{final} as providing a bound on the derivative of the POE (in $\beta>0$, when it is differentiable at $0$), by the $\widehat{\varphi}$-uniform expansion on  average condition. 
    \item The condition of weak unstable holonomies is used in the proof of Theorem \ref{final} in \eqref{strongMixLinBound}. Note that this assumption can be omitted when $F_\omega=F_{(\omega_i)_{i\leq0}}$.
\end{enumerate}
\end{remark}

\begin{cor}\label{rmkForStrongMixLinBoubd}
Under the assumptions of Theorem \ref{final}, for all $\omega_{\leq-1}\in \Sigma^-$,
$$\inf_{t\in X}\int \widehat{\varphi}\circ \widehat{F}^n(\omega,t) d\mu_{\omega_{\leq-1}}(\omega)\geq \chi -o_n(1),$$
and so also 
$$\inf_{t\in X}\int \frac{1}{n}\sum_{j\leq n}\widehat{\varphi}\circ \widehat{F}^j(\omega,t) d\mu_{\omega_{\leq-1}}(\omega)\geq \chi -o_n(1).$$
\end{cor}
This is a consequence of \eqref{strongMixLinBound}.

\appendix

\section{POE Maximizers and Future Independence}\label{appeB}

We saw in Theorem \ref{maxIsFutInd} that when $\wh\phi(\omega,t)=\phi(\omega)$, then any POE maximizer is future independent. This suggests a possible extension of the notion of stationary measures to Gibbs processes (see Remark \ref{newStationaryRmk}). Below in Proposition \ref{appAProp} we show that in the general case POE maximizers need not be future independent. In fact, we show an example where no POE maximizer is future independent.

\begin{prop}\label{appAProp}
 Let $\Sigma=\{0,1\}^\mathbb{Z}$. There exists a compact sub-shift of $\{0,1\}^\mathbb{Z}$, $(X,\sigma)$, 
and a potential $\widehat{\phi}\in \mathrm{H\ddot{o}l}(\Sigma\times X)$ with $\wh\phi(\omega,t)=\wh\phi(\omega_0,t_0)$, s.t. no POE maximizer is future independent, and in particular, $$\wh P(\wh \phi)< \sup\Big\{h_{\wh \eta \circ \pi^{-1}}(T)+\int \wh\phi d\wh\eta: \wh\eta\text{ is erg. }\wh F\text{-inv}\Big\}\leq P(\wh\phi),$$
 where $\pi:\Sigma\times X\to \Sigma$ is the projection, and $\wh F(\omega,t):=(T\omega, \sigma t)$ is simply the product dynamics. In this setting, $\wh\phi$ is trivially H\"older continuous, and trivially admits $\wh\phi$-holonomies.
\end{prop}
While the fiber dynamics and base dynamics are uncoupled, the potential (which is locally constant) does the coupling for the POE maximizers. Also, we in fact show that the maximizer in this example is unique.
\begin{proof}
By the Jewett-Krieger theorem (\cite{Jewett,Krieger}) there exists a uniquely ergodic sub-shift of $\{0,1\}^\mathbb{Z}$, $(X,\sigma,\tau)$, where $h_\tau(\sigma)>0$. Set $ \wh\phi(\omega,t)=C\cdot \mathbb{1}_{[\omega_0=t_0]}=C\cdot \Big(\mathbb{1}_{[1]\times [1]}+\mathbb{1}_{[0]\times [0]}\Big)$ where $C>0$ is a positive constant to be determined later.

We continue to bound the POE from above for all $\wh F$-invariant measures $\wh \eta$. Note that we are in the setting of Theorem \ref{VarPrincLowBound}, and so it is enough to bound the metric POE for every invariant measure. Note also that $\wh\eta$ must project to $\tau$ on $X$, as $(X,\sigma)$ is uniquely ergodic and our dynamics is a product dynamics. 

Let $\alpha:=\{\Sigma\times [0], \Sigma\times [1]\}$ and $N\in\mathbb{N}$, and write $\alpha_N=\bigvee_{i=-N}^N\wh F^{-i}\alpha$. Write $\wh{\mathcal{C}}_{-N}:=\bigvee_{i=-N}^0\wh F^{-i}\wh{\mathcal{C}}$. Then by the Martingales theorem (\cite{Martingales}), 

\begin{align}
    \wh P_{\wh\eta}(\wh \phi)=&\lim_{N\to\infty}\H_{\wh\eta}\Big(\wh{\mathcal{C}}|\wh{\mathcal{C}}_{-N}\vee \alpha_N\Big)+C\cdot (\wh\eta([1]\times [1])+\wh\eta([0]\times [0]))\label{appAeq0}\\
\leq&\H_{\wh\eta}\Big(\wh{\mathcal{C}}|\alpha\Big)+C\cdot (\wh\eta([1]\times [1])+\wh\eta([0]\times [0]))\label{appAeq1}\\
=&\sum_{b\in\{0,1\}}\tau([b])\cdot \Big(\H_{\wh\eta_b}\Big(\wh{\mathcal{C}}\Big)+C\frac{\wh\eta([b]\times [b])}{\wh\eta(\Sigma\times [b])}\Big)\nonumber\\
=&\sum_{b\in\{0,1\}}\tau([b])\cdot (\H(p_b)+C\cdot p_b)\nonumber\\
\leq& \H(p)+C\cdot p, \label{appAeq2}
\end{align}
where $\H(x):=-x\log x-(1-x)\log(1-x)$, $p_b:=\frac{\wh\eta([b]\times [b])}{\wh\eta(\Sigma\times [b])}$, $\wh\eta_b:=\frac{\wh\eta|_{\Sigma\times[b]}}{\tau([b])}$, and $p:=\frac{e^C}{1+e^C}\in(0,1)$ is the unique maximizer of the last equality. Also, note that, \eqref{appAeq1} is an equality if and only if $\wh{\mathcal{C}}$ is independent of $\bigvee_{i=-\infty}^{-1}\wh F^{-i}\wh{\mathcal{C}}\vee \bigvee_{i\neq0}\wh F^{-i}\alpha$ (w.r.t. $\wh\eta$), as the limit in \eqref{appAeq0} is non-increasing. Also, \eqref{appAeq2} is an equality if and only if $p_b=p$ for both $b=0,1$.

Next, we construct a maximizer $\wh\nu $ for the POE which obtains the upper bound in \eqref{appAeq2}. Given $t\in X$, we define a measure $\nu_t$ on $\Sigma\times\{t\}$ by,
\[
\nu_t \sim 
\begin{cases}
    \nu_t([\omega_i=t_i])=p  \\
    \nu_t([\omega_i=1-t_i])=1-p. 
\end{cases}
\]
We then define $\wh\nu:=\int \nu_td\tau(t)$. To see that $\wh\nu$ is indeed invariant, note that
$$\nu_t\circ \wh F^{-1}=\nu_{\sigma t}\Rightarrow \wh\nu\circ \wh F^{-1}=\wh\nu.$$
We now verify that $\wh\nu$ obtains the upper bound in \eqref{appAeq2}. First, for all $N\in\mathbb{N}$,
$$\H_{\wh\nu}\Big(\wh{\mathcal{C}}|\wh{\mathcal{C}}_{-N}\vee \alpha_N\Big)=\H_{\wh\nu}\Big(\wh{\mathcal{C}}|\alpha\Big),$$
as for every $t$, for $\nu_t$-a.e. $\omega$, $\omega_0$ is independent of $(\omega_i)_{i\leq-1}$ and of $(t_i)_{i\neq0}$. Therefore, $\wh P_{\wh\nu}(\wh \phi)=\H(p)+C\cdot p$.

It is left to show that $\wh\nu$ is the unique maximizer, and prove that it is not future independent. To show that it is the unique maximizer, assume that we have another maximizer $\wh\eta$. Then, equality must hold in both \eqref{appAeq1} and in \eqref{appAeq2}. That means, by the invariant of $\widehat{\eta}$, that for all $i\in\mathbb{Z}$, $\wh{\eta}(\omega_i=t_i|t_i)=p_{t_i,\wh\eta}=p$. Hence $\wh\eta=\wh\nu$. 

We now continue to show that $\wh\nu$ is not future independent. By Theorem \ref{boundedness}, it is equivalent to showing that $\wh h_{\wh\nu}(\wh F)<h_{\wh\nu\circ \pi^{-1}}(T)$. Indeed, by the fact that $\wh{\mathcal{C}}$ generates the $\sigma$-algebra of $\Sigma$, and $\alpha$ generates the $\sigma$-algebra of $X$, we have that 
\begin{equation}\label{appAeq3}
    h_{\wh\nu}(X|\Sigma)+h_{\wh\nu\circ\pi^{-1}}(T)=h_{\wh\nu}(\wh F)=h_{\wh\nu}(\Sigma|X)+h_\tau(\sigma),
\end{equation}
where $$h_{\wh\nu}(X|\Sigma)=\lim\limits_{N\to\infty}\frac{1}{N}\H_{\wh\nu}\Big(\alpha_N^+|\widehat{\mathcal{C}}_N\Big)\text{ and }h_{\wh\nu}(X|\Sigma)=\lim\limits_{N\to\infty}\frac{1}{N}\H_{\wh\nu}\Big(\widehat{\mathcal{C}}_N|\alpha_N^+\Big),$$ and $\alpha_N^+:=\bigvee_{i=0}^{N-1}\wh F^{-i}\alpha$.
First note that since $\omega_i$ depends only $t_i$, and so for all $N\in\mathbb{N}$, 
\begin{align*}
\H_{\wh\nu}\Big(\widehat{\mathcal{C}}_N|\alpha_N^+\Big)=&\H_{\wh\nu}\Big(\bigvee_{i=0}^{N-1}\wh F^{-i}\wh{\mathcal{C}}|\alpha_N^+\Big)=\H_{\wh\nu}\Big(\bigvee_{i=0}^{N-1}\wh F^{-i}\wh{\mathcal{C}}|\alpha_N^+\Big)\\
=&\H_{\wh\nu}\Big(\wh{\mathcal{C}}|\alpha_N^+\Big)+\H_{\wh\nu}\Big(\bigvee_{i=1}^{N-1}\wh F^{-i}\wh{\mathcal{C}}|\alpha_N^+\vee \wh{\mathcal{C}}\Big)\\
=&\H_{\wh\nu}\Big(\wh{\mathcal{C}}|\alpha_N^+\Big)+\H_{\wh\nu}\Big(\bigvee_{i=1}^{N-1}\wh F^{-i}\wh{\mathcal{C}}|\alpha_N^+\Big)\\
=&\cdots=\sum_{i=0}^{n-1}\H_{\wh\nu}\Big(\wh F^{-i}\widehat{\mathcal{C}}|\wh F^{-i}\alpha\Big)= N\cdot \H(p)=N\cdot \wh h_{\wh \nu}(\wh F).
\end{align*}

Plugging this back in \eqref{appAeq3}, we get that
\begin{equation}\label{appAeq4}
    h_{\wh\nu\circ\pi^{-1}}(T)=h_{\wh\nu}(\wh F)=h_{\wh \nu}(\wh F)+h_\tau(\sigma)-h_{\wh\nu}(X|\Sigma).
\end{equation}
Since $h_\tau(\sigma)>0$, if we can show that by choosing $C>0$ sufficiently large we can make $h_{\wh\nu}(X|\Sigma)$ arbitrarily small, we are done. We continue to bound $h_{\wh\nu}(X|\Sigma)$.

\begin{align}\label{appAeq5}
\H_{\wh\nu}\Big(\alpha_N^+|\widehat{\mathcal{C}}_N\Big)=&\H_{\wh\nu}\Big(\bigvee_{i=0}^{N-1}\wh F^{-i}\alpha|\widehat{\mathcal{C}}_N\Big)=\H_{\wh\nu}\Big(\bigvee_{i=0}^{N-1}\wh F^{-i}\alpha|\widehat{\mathcal{C}}_N\Big)\\
\leq&\sum_{i=0}^{N-1}\H_{\wh\nu}\Big(\wh F^{-i}\alpha|\widehat{\mathcal{C}}_N\Big)\leq \sum_{i=0}^{N-1}\H_{\wh\nu}\Big(\wh F^{-i}\alpha|\wh F^{-i}\widehat{\mathcal{C}}\Big)=N\cdot \H_{\wh\nu}\Big(\alpha|\widehat{\mathcal{C}}\Big).\nonumber
\end{align}

Set $\delta:=\wh\eta(\Sigma\times [1])=\tau([1])\in(0,1)$, and note that $\delta$ is independent of $C$, as every invariant measure project to $\tau$ on $X$ (by the unique ergodicity of $(X,\tau)$ and the product dynamics). A simple computation shows that for $\Phi(x):=-x\log x$,
\begin{align*}\label{appAeq6}
\H_{\wh\nu}\Big(\alpha|\widehat{\mathcal{C}}\Big)=&\Phi\Big(\frac{(1-\delta)(1-p)}{(1-\delta)(1-p)+(1-\delta)p}\Big)+\Phi\Big(\frac{\delta(1-p)}{\delta(1-p)+\delta p}\Big)\\
    +&\Phi\Big(\frac{(1-\delta)p}{(1-\delta)(1-p)+(1-\delta)p}\Big)+\Phi\Big(\frac{\delta p}{\delta(1-p)+\delta p}\Big).\nonumber
\end{align*}
It is clear that $p\xrightarrow[C\to\infty]{}1$, and so $\H_{\wh\nu}\Big(\alpha|\widehat{\mathcal{C}}\Big)\xrightarrow[C\to\infty]{}0$. This concludes the proof, by \eqref{appAeq4} and \eqref{appAeq5}. 

\end{proof}

\begin{remark}
    A slightly more detailed analysis could in fact guarantee also that, 
    \begin{equation}\label{appAeq7}
        \wh P(\wh\phi)<\sup\Big\{h_{\wh \eta \circ \pi^{-1}}(T)+\int \wh\phi d\wh\eta: \wh\eta\text{ is erg. }\wh F\text{-inv}\Big\}<P(\wh\phi).
    \end{equation}
    To see that, we notice that for every $\wh F$-invariant $\wh \eta$, for every $N\in\mathbb{N}$,
    \begin{align*}
\H_{\wh\eta}\Big(\alpha_N^+|\widehat{\mathcal{C}}_N\Big)=\H\Big(\alpha_N^+\vee\widehat{\mathcal{C}}_N\Big)-\H_{\wh\eta}\Big(\widehat{\mathcal{C}}_N\Big)\geq\H_{\wh\eta}\Big(\alpha_N^+\Big)-\H_{\wh\eta}\Big(\widehat{\mathcal{C}}_N\Big),
\end{align*}
and so, $h_{\wh\eta}(X|\Sigma)\geq h_\tau(\sigma)-h_{\wh\eta\circ\pi^{-1}}(T)$. Together with \eqref{appAeq3}, assuming that $h_\tau(\sigma)>\log 2$, we get that $h_{\wh\eta}(\wh F)\geq h_{\wh\eta\circ \pi^{-1}}(T)+(h_\tau(\sigma)-\log 2)$, and so \eqref{appAeq7} holds. However, to allow for $h_\tau(\sigma)>\log2$, we would have to choose $X$ as a sub-shift of $\{0,1,2\}^\mathbb{Z}$, which slightly complicates the already-complicated analysis Proposition \ref{appAProp}. 
\end{remark}

\section{Reduction to a one-sided shift}\label{appe}

\begin{lemma}\label{Sinai}
	Assume that $(\Sigma\times X,\wh F)$ admits stable holonomies. That is, for every $\omega,\omega'\in \Sigma$ which lie on the same stable leaf, the limit $\lim_{n\to\infty}(F_\omega^{n})^{-1}\circ F_{\omega'}^n\in \mathrm{H\ddot{o}l}_\alpha(X:X)$ converges in a H\"older norm with a uniform exponential rate. Then there exists a H\"older continuous mapping $\omega\mapsto F_\omega^+$ s.t. $F_{(\omega_i)_{i\in \mathbb{Z}}}^+= F_{(\omega_i)_{i\geq0}}^+ $ and $F_\omega^+$ is cohomologous to $F_\omega$, where the coboundary is $\alpha$-bi-H\"older continuous and depends in a H\"older manner on $\omega$.
\end{lemma}
\begin{proof}
For every $[a]\subset \Sigma$, fix $\omega_a\in [a]$. Denote by $[\omega_a,\omega]$ the Smale bracket. Define $C_\omega:= \lim_{n\to\infty}(F_\omega^{n})^{-1}\circ F_{[\omega_{\omega_0},\omega]}^n\in \mathrm{H\ddot{o}l}_\alpha(X:X)$. Note that $C_\omega$ is $\alpha$-bi-H\"older, as it admits an inverse $C_\omega^{-1}= \lim_{n\to\infty}F_{[\omega_{\omega_0},\omega]}^{-n}\circ F_\omega^{n} $. Then set $F_\omega^+:=C_{T\omega}^{-1}\circ F_\omega\circ C_\omega$, which is bi-H\"older.

We first show that $\omega\mapsto F_\omega^+$ is H\"older continuous. Let $\omega,\omega'\in \Sigma$ with $d(\omega,\omega)=e^{-N}<1$, and write $\omega_0=\omega_0'=a$. It is enough to show that $\omega\mapsto C_\omega$ is $\alpha$-H\"older continuous. By the uniform exponential convergence, let $K>0$ and and $\theta\in (0,1)$ s.t. 
\begin{align}\label{HoldOfCOmega}
	d_{\mathrm{H\ddot{o}l}_\alpha}(C_\omega, C_{\omega'})\leq& 2K e^{-\theta n}+ d_{\mathrm{H\ddot{o}l}_\alpha}((F_\omega^n)^{-1}\circ F^{n}_{[\omega_a,\omega]}, (F_{\omega'}^{n})^{-1}\circ F^{n}_{[\omega_a,\omega']}).
\end{align}
Choose $n=r\cdot N$, where $r\in (0,1)$ will be specified later. Let  $K,\tau>0$ s.t. $\omega\mapsto F_\omega$ is $(K',\tau)$-H\"older continuous. Then one can check by induction that $d_{\mathrm{H\ddot{o}l}_\alpha}((F_{\omega}^{n})^{-1}\circ F_{\omega'}^{n}, Id)\leq  B^ n\cdot  e^{-\tau (N-n)}$ for all $N$ large enough so $B^ n\cdot  e^{-\tau (N-n)} \leq 1$, where $B$ is a constant depending on\\ $\max_{\omega}\{\|F_\omega\|_{\mathrm{H\ddot{o}l}_\alpha},\|F_\omega^{-1}\|_{\mathrm{H\ddot{o}l}_\alpha}\}$, $K$, and $K'$. By choosing $r>0$ sufficiently small we can guarantee that $$d_{\mathrm{H\ddot{o}l}_\alpha}((F_{\omega}^{n})^{-1}\circ F_{\omega'}^{n}, Id)\leq K'' e^{-\frac{\tau}{2} N}.$$
Then, for all $N$ large enough (and $r>0$ sufficiently small),
 \begin{align*}
 	d_{\mathrm{H\ddot{o}l}_\alpha}((F_\omega^n)^{-1}\circ F^{n}_{[\omega_a,\omega]}, (F_{\omega'}^{n})^{-1}\circ F^{n}_{[\omega_a,\omega']})\leq K''' e^{-\frac{\tau}{3} N}.
 \end{align*}
Finally, putting this together with \eqref{HoldOfCOmega}, we conclude that $\omega\mapsto C_\omega$ is H\"older continuous.

We continue to show that $F_\omega^+$ depends only the non-negative coordinates of $\omega$:
	\begin{align*}
		F_\omega^+=&C_{T\omega}^{-1}\circ F_\omega\circ C_\omega\\
        =&\lim_n  F_{[\omega_{\omega_1},T\omega]}^{-n} \circ F_{T\omega}^{n}\circ F_\omega \circ (F_\omega^{n+1})^{-1}\circ F_{[\omega_{\omega_0},\omega]}^{n+1}\\
		=&\lim_n (F_{[\omega_{\omega_1},T\omega]}^{n})^{-1} \circ F_{[\omega_{\omega_0},\omega]}^{n+1},
	\end{align*}
	which depends only on $(\omega_i)_{i\geq0}$.
\end{proof}

The following lemma is proved similarly to Lemma \ref{Sinai}.

\begin{lemma}\label{Sinai2}
	Assume that $(\Sigma\times M,\wh F)$ admits stable holonomies, where $M$ is a closed Riemannian manifold, and $\wh F(\omega, t)=(T\omega, f_\omega(t))$ with $f_\omega\in \mathrm{Diff}^{1+\alpha}(M)$. That is, for every $\omega,\omega'\in \Sigma$ which lie on the same stable leaf, the limit $\lim_{n\to\infty}(f_\omega^{n})^{-1}\circ f_{\omega'}^n\in \mathrm{Diff}^{1+\alpha}(M)$ converges in $C^{1+\alpha}$ norm with a uniform exponential rate. Then there exists a H\"older continuous mapping $\omega\mapsto f_\omega^+$ s.t. $f_{(\omega_i)_{i\in \mathbb{Z}}}^+= f_{(\omega_i)_{i\geq0}}^+ $ and $f_\omega^+$ is conjugated to $f_\omega$, where the conjugacy is a $C^{1+\alpha}$ diffeomorphism which preserves fibers and it depends in a H\"older manner on $\omega$.
\end{lemma}

\bibliographystyle{alpha}
\tocless\bibliography{Elphi}

\def\cprime{$'$} \def\cprime{$'$}
\begin{thebibliography}{BEFRH}

\bibitem[BC14]{BergerCarrasco}
Pierre Berger and Pablo~D. Carrasco.
\newblock Non-uniformly hyperbolic diffeomorphisms derived from the standard map.
\newblock {\em Comm. Math. Phys.}, 329(1):239--262, 2014.

\bibitem[BEFRH]{BEFRH}
Aaron Brown, Alex Eskin, Simion Filip, and Federico Rodriguez-Hertz.
\newblock Measure rigidity for generalized u-gibbs states and stationary measures via the factorization method.
\newblock arXiv:2502.14042.

\bibitem[BG95]{BG_RandomRuelle}
Thomas Bogensch\"utz and Volker~Mathias Gundlach.
\newblock Ruelle's transfer operator for random subshifts of finite type.
\newblock {\em Ergodic Theory Dynam. Systems}, 15(3):413--447, 1995.

\bibitem[BKL23]{BuzziNonlinearTDF}
J\'er\^ome Buzzi, Beno\^it Kloeckner, and Renaud Leplaideur.
\newblock Nonlinear thermodynamical formalism.
\newblock {\em Ann. H. Lebesgue}, 6:1429--1477, 2023.

\bibitem[BOB25]{GuGibbs}
Snir Ben~Ovadia and David Burguet, 2025.

\bibitem[BQ16]{BenoistQuint}
Yves Benoist and Jean-Fran\c~cois Quint.
\newblock {\em Random walks on reductive groups}, volume~62 of {\em Ergebnisse der Mathematik und ihrer Grenzgebiete. 3. Folge. A Series of Modern Surveys in Mathematics [Results in Mathematics and Related Areas. 3rd Series. A Series of Modern Surveys in Mathematics]}.
\newblock Springer, Cham, 2016.

\bibitem[BRH17]{AaronFederico}
Aaron Brown and Federico Rodriguez~Hertz.
\newblock Measure rigidity for random dynamics on surfaces and related skew products.
\newblock {\em J. Amer. Math. Soc.}, 30(4):1055--1132, 2017.

\bibitem[CG12]{GouezelChazottes}
Jean-Ren\'e{} Chazottes and S\'ebastien Gou\"ezel.
\newblock Optimal concentration inequalities for dynamical systems.
\newblock {\em Comm. Math. Phys.}, 316(3):843--889, 2012.

\bibitem[DD]{DeWittDolgopyat}
Jon DeWitt and Dima Dolgopyat.
\newblock Expanding on average diffeomorphisms of surfaces: exponential mixing.
\newblock arXiv:2410.08445.

\bibitem[DH24]{YeorDima}
Dmitry Dolgopyat and Yeor Hafouta, 2024.

\bibitem[DK07]{DimaKrikorian}
Dmitry Dolgopyat and Rapha\"el Krikorian.
\newblock On simultaneous linearization of diffeomorphisms of the sphere.
\newblock {\em Duke Math. J.}, 136(3):475--505, 2007.

\bibitem[DKS08]{KiferDenkerStadlbauer}
Manfred Denker, Yuri Kifer, and Manuel Stadlbauer.
\newblock Thermodynamic formalism for random countable {M}arkov shifts.
\newblock {\em Discrete Contin. Dyn. Syst.}, 22(1-2):131--164, 2008.

\bibitem[ES23]{ELS}
Rosemary Elliott~Smith.
\newblock Uniformly expanding random walks on manifolds.
\newblock {\em Nonlinearity}, 36(11):5955--5972, 2023.

\bibitem[Jew70]{Jewett}
Robert~I. Jewett.
\newblock The prevalence of uniquely ergodic systems.
\newblock {\em J. Math. Mech.}, 19:717--729, 1969/70.

\bibitem[Kif01]{KiferRandomPressure}
Yuri Kifer.
\newblock On the topological pressure for random bundle transformations.
\newblock In {\em Topology, ergodic theory, real algebraic geometry}, volume 202 of {\em Amer. Math. Soc. Transl. Ser. 2}, pages 197--214. Amer. Math. Soc., Providence, RI, 2001.

\bibitem[Kri72]{Krieger}
Wolfgang Krieger.
\newblock On unique ergodicity.
\newblock In {\em Proceedings of the {S}ixth {B}erkeley {S}ymposium on {M}athematical {S}tatistics and {P}robability ({U}niv. {C}alifornia, {B}erkeley, {C}alif., 1970/1971), {V}ol. {II}: {P}robability theory}, pages 327--346. Univ. California Press, Berkeley, CA, 1972.

\bibitem[LW77]{LedWalt}
Fran\c~cois Ledrappier and Peter Walters.
\newblock A relativised variational principle for continuous transformations.
\newblock {\em J. London Math. Soc. (2)}, 16(3):568--576, 1977.

\bibitem[Pot22]{Rafael}
Rafael Potrie.
\newblock A remark on uniform expansion.
\newblock {\em Rev. Un. Mat. Argentina}, 64(1):11--21, 2022.

\bibitem[PS82]{uGibbsPesinSinai}
Ya.~B. Pesin and Ya.~G. Sina\u{\i}.
\newblock Gibbs measures for partially hyperbolic attractors.
\newblock {\em Ergodic Theory Dynam. Systems}, 2(3-4):417--438 (1983), 1982.

\bibitem[Rue67]{Ruelle67}
D.~Ruelle.
\newblock A variational formulation of equilibrium statistical mechanics and the {G}ibbs phase rule.
\newblock {\em Comm. Math. Phys.}, 5:324--329, 1967.

\bibitem[Wal75]{Walters}
Peter Walters.
\newblock A variational principle for the pressure of continuous transformations.
\newblock {\em Amer. J. Math.}, 97(4):937--971, 1975.

\bibitem[Wil91]{Martingales}
David Williams.
\newblock {\em Probability with Martingales}.
\newblock Cambridge University Press, 1991.

\end{thebibliography}

\Addresses

\end{document}